
\input amstex
\documentstyle{amsppt}
\magnification 1200
\def\boxit#1{\vbox{\hrule\hbox{\vrule\kern2.5pt\vbox{\kern2.5pt#1
\kern2.5pt}\kern2.5pt\vrule}\hrule}}
\NoBlackBoxes
\input epsf.tex
\topmatter
\title Higher Syzygies of  Elliptic Ruled Surfaces \endtitle
\author Francisco Javier Gallego $^*$ \\ and \\ B. P. Purnaprajna
\endauthor
\address{F.J. Gallego, Dpto. de Algebra, Facultad de Matem\'aticas, U.C.M.,
28040 Madrid SPAIN}\endaddress
\email gallego\@sunal1.mat.ucm.es\endemail
\address{B.P. Purnaprajna, Dept. of Mathematics, Brandeis
University, Waltham MA 02254-9110 USA.}\endaddress
\email purna\@max.math.brandeis.edu\endemail
\date{10 June 1995}\enddate
\thanks{We would like to thank our advisor David Eisenbud for
his  help, patience and encouragement. We would also like to
thank Robert Lazarsfeld and Mohan Kumar for their encouragement 
and advice. $*$
Partially supported by DGICYT,
PB93-0440-C03-01.}\endthanks
\endtopmatter
\document
\headline={\ifodd\pageno\rightheadline \else\leftheadline\fi}
\def\leftheadline{\tenrm\hfil {\eightpoint F.J.GALLEGO \& B.P. PURNAPRAJNA}
\hfil\folio}
\def\rightheadline{\tenrm\folio\hfil {\eightpoint HIGHER SYZYGIES OF ELLIPTIC
RULED SURFACES}
\hfil}
\voffset=2\baselineskip
\vsize= 22truecm
\nologo
\vskip .3cm

\roster
\item ""{\bf  Introduction}
\item "{\bf 1.}"{\bf Background material}
\item "{\bf 2.}" {\bf An example and a general result}
\item "{\bf 3.}" {\bf Some  lemmas and
commutative diagrams}
\item "{\bf 4.}" {\bf Cohomology vanishings on ruled
elliptic surfaces with invariant $e = -1$}
\item"{\bf 5.}" {\bf Cohomology vanishings on
elliptic ruled surfaces with invariant $e \geq
0$}
\item "{\bf 6.}" {\bf Syzygies of elliptic ruled
surfaces}
\item"{\bf 7.}" {\bf Open questions and conjectures}
\endroster

\heading   Introduction\endheading

The purpose of this article is to study the minimal free resolution of
homogeneous
coordinate rings of elliptic ruled surfaces.

 Let $X$ be an irreducible projective
variety and $L$ a very ample line bundle on $X$, whose
complete linear series defines the morphism
$$\phi_L : X \longrightarrow \bold P (\text H^0(L))$$
 Let $S= \bigoplus_{m=0}^\infty\text S^m\text H^0(X,L)$ and $R(L)\bigoplus_{m=0}^\infty\text H^0(X, L^{\otimes m})$. Since $R(L)$ is
a finitely generated  graded module over S,  it has a minimal
graded free resolution. We say that the line bundle $L$ is normally
generated if the natural maps
$$ \text S^m\text H^0(X,L) \to \text H^0(X,L^{\otimes m}) $$
are surjective for all $m \geq 2$. If $L$ is normally generated,
then we say $L$ satisfies property $N_p$,
if the matrices in the free resolution of $R$ over $S$ have linear
entries until the pth stage.

In this article we prove the following result (Theorem 6.1): Let $X$
be an elliptic ruled surface  and let $L=B_1\otimes ...\otimes
B_{p+1}$ be a line bundle on $X$, such that each $ B_i$ is
base-point-free and ample. Then $L$ satisfies property $N_p$. As a
corollary of this result we show that the adjoint bundle
$\omega _X
\otimes A_1 \otimes \dots \otimes A_{2p+3}$ satisfies property
$N_p$, for arbitrary ample line bundles $A_i$.

To put things in perspective, we would like to recall some well
known results in this
area. On the subject of adjoint linear series,   Reider recently
proved (c.f. [R]) that if
$X$ is a surface over the complex numbers and $A$ is an ample line
bundle, then
$\omega_X
\otimes A ^{\otimes 4}$ is very ample. Mukai has conjectured that
$\omega_X
\otimes A ^{\otimes p + 4}$ satisfies  property $N_p$.
 Some work in this direction
has been done by David Butler in [Bu], where he studies the syzygies
of adjoint linear series  on ruled varieties. He proves
that if the dimension of $X$ is $n$, then $\omega _X
\otimes A ^{\otimes 2n+1}$ is normally generated and $\omega _X
\otimes A ^{\otimes 2n+2np}$ satisfies property $N_p$.
In particular if $X$ is a ruled surface
Butler's result says that  $\omega _X
\otimes A ^{\otimes 5}$ is normally generated and $\omega _X
\otimes A ^{\otimes 4+4p}$ satisfies property $N_p$.

Of course a more general
question would be: given a very ample line bundle $L$ on $X$, what
is the largest
$p$ such that
 property $N_p$ holds for $L$?
 A relevant result in this line is due to Yuko
Homma (c.f. [Ho1] and [Ho2]),
 who has classified  all line bundles which are normally generated
on an elliptic ruled surface. Another result is obtained in [GP],
where we characterize those line bundles which satisfy property
$N_1$ on an elliptic ruled surface. Therefore, it is in the light
of the above more general question that Theorem 6.1 should be
regarded. Certainly the mentioned results provide information about
the particular case of adjoint linear series. The result of Homma
and our results in [GP] imply Mukai's conjecture, in the case of
elliptic ruled surfaces, for $p=0$ and $p=1$ respectively.
Accordingly, as a
corollary of Theorem 6.1 we improve the bound obtained by Butler by
almost a factor of two in the case of  elliptic ruled
surfaces.

One of the tools we use in  this article is the so-called Koszul
cohomology, developed by Mark Green, which links the study of the
vanishing of graded Betti numbers of minimal free resolutions to the
study of cohomology  vanishings of certain vector bundles.
We also use a theorem of Castelnuovo
on surjectivity of multiplication maps of global sections,
generalized by Mumford, (c.f. Theorem 1.3). Then in Sections 2 and 3
we develop the machinery that will allow us to prove our
main results.
For the set-up of this machinery we
use induction on the number of base-point-free divisors and on the
dimension of the variety. Briefly, we will associate to a line
bundle, which is a product of certain number of base-point-free
line bundles, another bundle. We will show that the  first
cohomology group of the associated bundle vanishes. To achieve this
goal for a variety
$X$ of arbitrary dimension we restrict the line bundle to a
divisor of $X$ and then use induction on the dimension. And to
achieve the goal for an arbitrary number of base-point-free divisors
in the mentioned product, we use induction on the number of
base-point-free divisors.

The usefulness of these constructions is
not limited to the example of elliptic ruled surfaces. In this article
(see
Theorem 2.2) we also obtain from them, results for all surfaces with
geometric genus
$0$ (a class that includes  Enriques surfaces).
In fact our results can be summarized in the following
principle: let $L$ be the tensor product of $p+1$ ample,
base-point-free line bundles; If certain cohomology vanishings occur
then $L$ should satisfy property $N_p$. This principle holds in
wider generality: In two forthcoming articles ([GP1],[GP2]) we use the
machinery
developed here to show that the  principle mentioned above holds for
surfaces with
$\kappa =0$, Fano varieties of dimension $n$ with index bigger than
or equal to
$n-2$ and elliptic surfaces. We also show that the principle holds,
with some minor modifications which depend on the dimension of the
variety, for pluricanonical models of surfaces of general type and
Calabi-Yau threefolds. The former answers an open question in [B].

\heading  1. Background material \endheading
\proclaim {Convention} Throughout this paper we will work over an
algebraically closed \text{field $\bold k$.}
\endproclaim

One of the tools we will use is  this beautiful cohomological characterization
by Green of the property $N_p$
(c.f. [G], [GL], [L]).
Let
$L$ be a globally generated line bundle. We define the vector bundle $M_L$ as
follows:
$$0 \to M_L \to \text H ^0(L) \otimes \Cal O \to
L \to 0\eqno (1.1)$$
\proclaim {Lemma 1.2} Let $L$ be a normally generated,
nonspecial line bundle on a surface $X$ with geometric genus
$0$. Then, $L$ satisfies the property $N_p$ iff $\text H ^1(\wedge
^{p'+1} M_L \otimes L) = 0$ for all $1 \leq p' \leq p$.
\endproclaim
{\it Proof.} See [GL] \S 1.
\boxit{}

(1.2.1) If the characteristic of $\bold k$ is  strictly bigger than
$p+1$, then  the vanishings of  $\text H ^1(\wedge ^{\otimes  p'
+1} M_L\otimes L)$ for all $0 \leq p' \leq p$ will follow from  the
vanishings of
$\text H ^1(M_L^{\otimes  p' +1}\otimes L)$, because $\wedge ^{\otimes  p'
+1} M_L\otimes L$ is a direct summand of $M_L^{\otimes  p'
+1}\otimes L$.

\vskip .2 cm

The other main tool we will use is a generalization by Mumford
of a lemma of Castelnuovo:
\proclaim {Theorem  1.3} Let $L$ be a base-point-free line bundle on a variety
$X$ and let $\Cal F$ be a coherent sheaf on $X$. If $\text H^i (\Cal F \otimes
L ^{-i}) = 0$, for all $i \geq 1$, then the multiplication map
$$ \text H ^0(\Cal F \otimes L ^{\otimes i}) \ \otimes \text
H ^0(L) \to  \text H ^0(\Cal F \otimes L ^{\otimes i+1})$$
is surjective, for all $i \geq 0$.
\endproclaim

{\it Proof.} [Mu], p. 41, Theorem 2. Note that the assumption there
of $L$ being ample is unnecessary. \boxit{}

\vskip .2 cm
It will be useful to have the following characterization of projective
normality:
\proclaim {Lemma 1.4} Let $X$ be a surface with geometric genus $0$ and let
$L$ be an ample, base-point-free line bundle.
If $\text H^1(L )
=0$, then $L$ is normally generated iff $\text H^1(M_L \otimes
L) =0$.
\endproclaim

{\it Proof.} See [GP], Lemma 1.4. \boxit{}

\heading  2. An example and a general result \endheading

In [GP] we gave a complete characterization of the divisors on
an elliptic ruled surface satisfying the property
$N_1$. In particular we showed ([GP], Corollary 4.4) that
$$B ^{\otimes 2} \ \text{ satisfies the  property} \ N_1
\ \text{if $B$ is  ample and base-point-free.} \leqno
(2.1)$$   To give an idea of how we will generalize to higher
syzygies the results and techniques from [GP] we will focus in this
section on the generalization of (2.1). We recall that the statement
in  (2.1)  was shown to be true for a larger class of surfaces (see
[GP], Proposition 2.1 and Corollary 2.8), namely, those with $p_g = 0$,
if one requires $B$ to be nonspecial (the latter condition is
automatically satisfied by ample base-point-free line bundles on
elliptic ruled surfaces). Thus, we will prove the following

\proclaim{Theorem 2.2} Let
$X$ be  a surface with
$p_g = 0$. Let $\text{\rm{char}} \ \bold k > p+1$ or equal to $0$.
Let
$B$ be a nonspecial, ample, and base-point-free line bundle. Then
$B^{\otimes p +1}$ satisfies the property $N_p$, for all $p \geq
1$.\endproclaim

(2.2.1) The  same statement is false for $p=0$. Consider for instance $X$
elliptic ruled surface of invariant $e=-1$, let $C_0$ be a minimal section
of $X$ and let
$B$  equal
$\Cal O_X(2C_0)$. $B$ is an ample, base-point-free line bundle (c.f.
Propositions 3.4 and 3.5) but it is not very ample since its
restriction to
$C_0$ is not very ample.

Before we prove Theorem 2.2 we will require the
following

\proclaim{Lemma 2.3} Let $X$ and $B$ be as in Theorem
2.2. If $p \geq 1$ and $p_1$,$p_2\geq
p$, then the cohomology groups
\text{$\text H ^1(M_{B ^{\otimes
p_1+1}}^{\otimes p+1}
\otimes B ^{\otimes p_2+1})$} and $\text H ^2(M_{B
^{\otimes p_1+1}}^{\otimes p+1}
\otimes B ^{\otimes p_2})$ vanish.\endproclaim

Before we give the proof of Lemma 2.3 we make two observations
\proclaim{Observation 2.4}
Let $X $ be a surface with geometric genus 0,
let $B$ be a base point free line bundle and let $P$ be an
effective line bundle such that H$^1(P)= $H$^1(B) =0$. Then
H$^1(B\otimes P)=0$.
\endproclaim

\proclaim
{Observation 2.5} Let $X$ be a surface, let $P$ be an
effective line bundle and L a coherent sheaf. If H$^2(L) =0$, then
H$^2(L \otimes P) =0$
\endproclaim

{\it (2.6) Proof of Lemma 2.3.}
The proof is by induction on $p$. If $p =1$ we have to prove
that
$$\matrix (2.6.1) && \text H^1(M_{B^{\otimes
a}}^{\otimes 2}
\otimes B^{\otimes b}) &=& 0 \cr
(2.6.2) && \text H^2(M_{B^{\otimes
a}}^{\otimes 2}
\otimes B^{\otimes b -1}) &=& 0 \cr
\endmatrix$$
for all $a$,$b \geq 2$.
Note that  if $\text H ^1(M_{B ^{ \otimes a}} \otimes B
^{\otimes b}) = 0$, the vanishing in  (2.6.1) is equivalent to
the surjectivity of the following multiplication map:
$$\text H ^0(M_{B ^{ \otimes a}} \otimes B
^{\otimes b}) \otimes \text H ^0(B^{\otimes a}) @>\alpha>>
\text H ^0(M_{B ^{ \otimes a}} \otimes B
^{\otimes a+b}) \ .$$
To show the surjectivity of $\alpha$ it suffices to show the
surjectivity of
$$\text H ^0(M_{B ^{ \otimes a}} \otimes B
^{\otimes b}) \otimes \text H ^0(B)^{\otimes a} \to
\text H ^0(M_{B ^{ \otimes a}} \otimes B
^{\otimes a+b}) \ .$$ From all the above and from Theorem 1.3, it
follows that in order to prove (2.6.1),  it is enough to show
that
$$\displaylines{\rlap{(2.6.3)} \hfill\text H ^1(M_{B ^{ \otimes a}}
\otimes B
^{\otimes b'-1}) = 0 \hfill\cr
 \rlap{(2.6.4)}\hfill \text H ^2(M_{B ^{ \otimes a}} \otimes B
^{\otimes b'-2}) = 0 \hfill\cr}$$ for all $a$,$b' \geq 2$. From
Observation 2.4 it follows that $\text H ^1(B ^{\otimes c}) = 0 $
for all $c
\geq 1$. Therefore the  vanishing  (2.6.3) is equivalent to the surjectivity
of the map
$$\text H ^0(B^{\otimes a}) \otimes \text H ^0( B
^{\otimes b'-1})   \to
\text H ^0( B
^{\otimes a+b'-1}) \ .\eqno (2.6.5)$$ Thus it suffices to show
the surjectivity of
$$\text H ^0(B^{\otimes r}) \otimes \text H ^0( B)
^{\otimes s}   \to
\text H ^0( B
^{\otimes r+s})$$
for all $r$ and $s$ such that $r \geq s$, $r \geq 2$ and $s \geq
1$. This follows at once from Theorem 1.3, since, by Observations
2.4 and 2.5, $\text H ^1(B^{\otimes c}) = \text H ^2(B^{\otimes d})
=0$
for all
$c
\geq 1$ and all $d \geq 0$.

{}From the exact sequence 1.1, it follows that the vanishings of
both
$\text H ^1(B^{\otimes a +b'-2})$ and $\text H ^2(B^{\otimes b'-2})
$ imply (2.6.4). Now we
prove (2.6.2). From the exact sequence 1.1 it is enough to show
that
$\text H ^1(M_{B ^{ \otimes a}}
\otimes B
^{\otimes a+b-1})$
and $\text H ^2(M_{B ^{ \otimes a}} \otimes B
^{\otimes b-1})$ vanish. These vanishings are special cases of
(2.6.3) and (2.6.4).

Now assume that the result is true for $p-1$.
Then in particular \linebreak \text{$\text
H ^1(M_{B ^{\otimes p_1+1}}^{\otimes
p}
\otimes B ^{\otimes p_2+1}) = 0$} and therefore the vanishing of
the group
\linebreak \text{$\text H ^1(M_{B ^{\otimes p_1+1}}^{\otimes p+1}
\otimes B ^{\otimes p_2+1}) $} is equivalent to the
surjectivity of the following map:
$$\text H ^0(M_{B ^{\otimes p_1+1}}^{\otimes p}
\otimes B ^{\otimes p_2+1}) \otimes \text H ^0(B^{\otimes p_1+1})
@>\gamma>>
\text H ^0(M_{B ^{\otimes p_1+1}}^{\otimes p}
\otimes B ^{\otimes p_1 +p_2+2}) \ .$$
The surjectivity of $\gamma$ follows from the surjectivity of
$$\text H ^0(M_{B ^{\otimes p_1+1}}^{\otimes p}
\otimes B ^{\otimes p_2+1}) \otimes \text H ^0(B)^{\otimes p_1+1}
\to
\text H ^0(M_{B ^{\otimes p_1+1}}^{\otimes p}
\otimes B ^{\otimes p_1 +p_2+2})$$ and this in turn follows
from Theorem 1.3 and induction hypothesis.

To show that $\text H
^2(M_{B ^{\otimes p_1+1}}^{\otimes p+1}
\otimes B ^{\otimes p_2}) = 0$ it suffices, again by the exact
sequence 1.1, to
check that the groups $\text H ^1(M_{B
^{\otimes p_1+1}}^{\otimes p}
\otimes B ^{\otimes p_1+p_2+1})$ and $\text H ^2(M_{B
^{\otimes p_1+1}}^{\otimes p}
\otimes B ^{\otimes p_2})$ vanish. This is true by induction.
\boxit{}

\vskip .45 cm

(2.7) {\it Proof of Theorem 2.2.}
By Lemmas 1.2 and 1.4 and by (1.2.1) it is enough to show that
$\text H ^1(M_{B ^{\otimes p+1}}^{\otimes  p' +1}\otimes B
^{\otimes p+1})$ vanishes for all
$0
\leq p'
\leq p$. The vanishing follows when $1 \leq p' \leq p$ as a
particular case of Lemma 2.3. Since $\text H ^1 ( B ^{\otimes c})
= 0$ for all $c \geq 1$, the vanishing of $\text H
^1(M_{B ^{\otimes p+1}}\otimes B ^{\otimes p+1})$ is equivalent
to the surjectivity of the multiplication map
$$\text H ^0(B^{\otimes p+1}) \otimes \text H ^0( B
^{\otimes p+1})   \to
\text H ^0( B
^{\otimes 2p+2}) \ .$$
and that is a special case of (2.6.5). \boxit{}
\proclaim{Corollary 2.7.1}
Let $X$ be an Enriques surface over an algebraic closed field of
characteristic $0$. Let
$B$ be an ample base-point-free line bundle. Then $B^{\otimes p+1}$
satisfies the property $N_p$, for all $p \geq 1$.
\endproclaim
{\it Proof.}
Since $K_X \equiv 0$ and $B$ is ample, $\omega_X \otimes B$ is also
ample and by Kodaira vanishing, H$^1(B)=0$. Thus we can apply
Theorem 2.2. \boxit{}

\vskip .3 cm

In Theorem 2.2 we have dealt with line bundles which are powers of
a base-point-free line bundle. Obviously not all the line bundles
on a surface $X$ are of this form. Therefore we want now to study the
syzygies of a wider variety of line bundles. For this
purpose it is convenient to abstract and  somehow generalize the
formalism of Lemma 2.3. We will do so in the next lemma, which is
key to the proof of Propositions 4.1, 4.2 and 5.1, on which the
results of Section 6 are based.

\proclaim{Lemma 2.8} Let $X$ be a surface. Let $q_0$ be
a positive integer and let $\frak B$ and $\frak P$ be two subsets of
$\text{Pic}(X)$ satisfying the following properties:
\roster
\item  "2.8.1." All elements in $\frak B$ are base-point-free and
if $B \in
\frak B$ and $B \equiv B'$, then $B' \in \frak B$. The set $\frak B
$ is contained in $\frak P$ and if $P^1$ and $P^2$ belong to $\frak
P$, then
$P^1 \otimes P^2$ belongs to $\frak P$.

\item  "2.8.2."  For all $B_1, \dots , B_{q_0+3} \in \frak B$, the
line bundle $B_1 ^{\otimes 2} \otimes B_2 \otimes \dots \otimes
B_{q_0+2} \otimes B_{q_0+3}^*$ belongs to $\frak P$.

\item "2.8.3." For all
$B_1, \dots ,
B_{q_0+1}, B_1', \dots , B_{q_0+2}', C_1, \dots, C_n \in \frak B$
such that
$B_i \equiv B_i'$ and for any
 line bundle $P \in \frak P$, the line bundles $$\displaylines{ R_3 = B_1
\otimes
\dots
\otimes B_{q_0+1} \otimes C_1 \otimes \dots \otimes C_n \quad
\text{and}\cr R_3' = B_1' \otimes \dots
\otimes B_{q_0+1}' \otimes P \otimes B^{\prime *}_{q_0+2}\cr}$$
satisfy
$\text H ^2(M_{R_3}^{\otimes q_0+1}\otimes R_3') = 0$.

\item "2.8.4." For all  $B_1, \dots ,
B_{q_0+1}, B_1', \dots ,
B_{q_0+1}', C_1, \dots, C_n \in \frak B$ such that
$B_i \equiv B_i'$ and for any
line bundle $P \in \frak P$, the line bundles $$\displaylines{R_4 = B_1 \otimes
\dots
\otimes B_{q_0+1} \otimes C_1 \otimes \dots \otimes C_n
\ \text{and}\cr  R_4' = B_1' \otimes \dots \otimes B'_{q_0+1}
\otimes P\cr}$$ satisfy
$\text H ^1(M_{R_4}^{\otimes q_0+1}\otimes R_4') = 0$.
\endroster
Given $q \geq q_0$, let \ $T_1, \dots , T_{q+1}, T_1', \dots ,
T_{q+1}', S_1,
\dots , S_n
\in \frak B$ such that $T_i \equiv T_i'$ and let $Q \in \frak P$. If $R = T_1
\otimes \dots \otimes T_{q+1} \otimes S_1 \otimes \dots \otimes
S_n$ and $R' = T_1' \otimes \dots \otimes T_{q+1}' \otimes Q$, then
$$ \text H^1(M_R^{\otimes q+1}\otimes R')=0 \ .$$
\endproclaim

{\it Proof.}
We prove the lemma using induction on $q$. If $q =q_0$ the result is
just Condition 2.8.4. Now assume that the result is true for $q_0,
\dots , q-1$.  After tensoring  exact sequence 1.1 by $M_R^{\otimes
q} \otimes R'$ and taking global sections we obtain:
$$\displaylines {\text H^0(M_R ^{\otimes q} \otimes R') \otimes \text
H^0(R) @>\alpha>> \text H^0(M_R ^{\otimes q} \otimes R' \otimes R) \cr
\to
\text H^1(M_R ^{\otimes q+1} \otimes R') \to \text H^1(M_R ^{\otimes q} \otimes
R') \otimes \text H^0(R)}$$
Using 2.8.1 and induction hypothesis on $q-1$ it follows that
\text{$\text H^1(M_R
^{\otimes q} \otimes R')$} vanishes. Thus the surjectivity of
$\alpha$ is equivalent to the vanishing of the group $\text H^1(M_R
^{\otimes q+1}
\otimes R')$. We argue like this: The surjectivity of
$\alpha$ follows from the surjectivity of
$$\text H^0(M_R ^{\otimes q} \otimes R') \otimes
\bigotimes_{i=1}^{q+1}\text H^0(T_i)
\otimes \bigotimes_{j=1}^n
\text H^0(S_j)
 @>\beta >> \text H^0(M_R ^{\otimes q} \otimes
R'
\otimes R)
$$
and to obtain the surjectivity of $\beta$, by Theorem 1.3,
it is sufficient to check the following vanishings:
$$\displaylines{\rlap{(2.8.5)}\hfill\hfill\cr
\hfill
\text H^1(M_R^{\otimes q} \otimes T_1^{\otimes 2} \otimes \dots
\otimes T_{i-1}^{\otimes 2} \otimes T_{i+1} \otimes \dots \otimes
T_{q+1}
\otimes N \otimes Q)=0\hfill\cr
\hfill\text{for all}\ 1 \leq i \leq q+1 \ \text{and any} \ N \
\text{nef}\hfill\cr
\rlap{(2.8{.6})}\hfill\hfill\cr
\hfill\text H^2(M_R^{\otimes q} \otimes
T_1^{\otimes 2}
\otimes
\dots
\otimes T_{i-1}^{\otimes 2} \otimes T_{i+1} \otimes \dots \otimes
T_{q+1}
\otimes T_i^*
\otimes N \otimes Q)=0\hfill\cr
\hfill\text{for all} \ 1 \leq i \leq q+1 \ \text {and any}\ N \
\text{nef}\hfill\cr
\rlap{(2.8.7)}\hfill\hfill\cr
\hfill\text H^1(M_R^{\otimes q} \otimes
T_1^{\otimes 2} \otimes
\dots
\otimes  T_{q+1}^{\otimes 2} \otimes S_1 \otimes \dots \otimes
S_{j-1}
\otimes S_j^*
\otimes N \otimes Q)=0\hfill\cr
\hfill\text{for all}\ 1 \leq j \leq
n \
\text {and any} \
N \ \text{nef}\hfill\cr
\rlap{(2.8.8)}\hfill\hfill\cr
\hfill\text H^2(M_R^{\otimes q} \otimes
T_1^{\otimes 2} \otimes
\dots
\otimes  T_{q+1}^{\otimes 2} \otimes S_1 \otimes \dots \otimes
S_{j-1}
\otimes S_j^{-2}
\otimes N \otimes Q)=0 \hfill\cr
\hfill\text{for all} \ 1 \leq j \leq
n \
\text {and any} \
N \ \text{nef}
\hfill\cr}$$
First we check (2.8{.6}). If $q=q_0+1$, (2.8{.6}) follows from
Conditions 2.8.1 ($B \otimes N \in \frak B$ for all $B \in \frak
B$) and 2.8.3. If $q \geq q_0+2$, using exact sequence 1.1 it
suffices to check that
$$\displaylines{\rlap{(2.8.9)}\hfill\hfill\cr
\hfill
\text H^1(M_R^{\otimes k} \otimes R \otimes T_1^{\otimes 2} \otimes
\dots
\otimes T_{i-1}^{\otimes 2} \otimes T_{i+1} \otimes \dots \otimes
T_{q+1}
\otimes T_i^* \otimes N \otimes Q)=0 \hfill\cr
\hfill\text{for all}\
1
\leq i
\leq q+1 ,\
\text{for all}\  q_0+1
\leq k \leq q-1 \ \text {and any} \
N \ \text{nef}\hfill\cr
\rlap{(2.8.10)}\hfill\hfill\cr
\hfill\text H^2(M_R^{\otimes q_0+1} \otimes
T_1^{\otimes 2}
\otimes
\dots
\otimes T_{i-1}^{\otimes 2} \otimes T_{i+1} \otimes \dots \otimes
T_{q+1}
\otimes T_i^*
\otimes N \otimes Q)=0\hfill\cr
\hfill\text{for all} \ 1 \leq i
\leq q+1 \ \text {and any}\
N \ \text{nef} \hfill\cr}$$
The vanishing in (2.8.10) follows from Condition 2.8.1 and Condition 2.8.3.
We will postpone the proof of (2.8.9) for the moment. Now we check
(2.8.8):
If $q=q_0+1$, then (2.8.8) follows from Conditions 2.8.1, 2.8.2 and
2.8.3. If
$q \geq q_0+2$, again using exact sequence 1.1 it suffices to
check that
$$\displaylines{\rlap{(2.8.11)}\hfill\hfill\cr
\hfill\text H^1(M_R^{\otimes k} \otimes R \otimes T_1^{\otimes 2}
\otimes
\dots
\otimes  T_{q+1}^{\otimes 2} \otimes S_1 \otimes \dots \otimes
S_{j-1}
\otimes S_j^{-2}
\otimes N \otimes Q)=0 \hfill\cr
\hfill\text{for all} \ 1 \leq j \leq n ,\
\text{for all} \  q_0+1
\leq k \leq q-1 \ \text {and any}\
N \  \text{nef}\hfill\cr
\rlap{(2.8.12)}\hfill\hfill\cr
\hfill \text H^2(M_R^{\otimes q_0+1}
\otimes T_1^{\otimes 2}
\otimes
\dots
\otimes  T_{q+1}^{\otimes 2} \otimes S_1 \otimes \dots \otimes
S_{j-1}
\otimes S_j^{-2}
\otimes N \otimes Q)=0\hfill\cr
\hfill\text{for all} \ 1 \leq j \leq n\
\text {and any}\
N \ \text{nef}\hfill\cr}
$$
The vanishing in (2.8.12) follows from Conditions 2.8.1, 2.8.2, and 2.8.3.

We still have to check (2.8.5), (2.8.7), (2.8.9) and (2.8.11). The vanishing
in (2.8.5) follows from Condition 2.8.1 and induction hypothesis on $q-1$.
The vanishing in (2.8.7) follows from Conditions 2.8.1 and 2.8.2 and
induction hypothesis on $q-1$. The vanishings in (2.8.9) and (2.8.11) follow
from Conditions 2.8.1 and 2.8.2 and induction hypothesis on $q_0,
\dots , q-2$.
\boxit{}

\heading  3. Some  lemmas and
commutative diagrams \endheading
In this section we prove several lemmas which we will   use in
Sections 4 and 5. The first three lemmas hold in great generality.
The first one is connected to this problem: Consider
two base-point-free
line bundles $L_1$ and $L_2$. We would like to relate  the
vanishing of
the cohomology of $M_{L_1}^{\otimes p+1}
\otimes L_2$ to the vanishing of the cohomology of a similar bundle
on a divisor
$Y$ of
$X$, obtained by restricting $L_1$ and $L_2$ to $Y$. The second
and third lemma deal roughly with the following situation:
 Consider in addition to $L_1$ and $L_2$,
  two ``bigger" line bundles $L_1'$ and $L_2'$  (in the sense that
\text{$L_i' \otimes L_i^*$} is an effective line bundle). We would
like to
 relate the vanishing of the cohomology of $M_{L_1'}^{\otimes p+1}
\otimes L_2$ and $M_{L_1}^{\otimes p+1}
\otimes L_2'$ to the vanishing of the cohomology of
$M_{L_1}^{\otimes p+1}
\otimes L_2$. The usefulness of these kinds
of results is quite clear. For example, they give us a way to prove
 that if a line bundle $L$ satisfies the property $N_p$, then
so does the tensor product of $L$ with certain effective line
bundles. Therefore these three lemmas will be  a key
element in the proofs of Propositions 4.1, 4.2 and 5.1.

\proclaim{Lemma 3.1} Let $X$ be a projective variety, let $q$ be a
nonnegative integer and let $F_i$ be a base-point-free line bundle
on $X$ for all $1 \leq i \leq q+1$. Let $Q$ be an effective line
bundle on $X$ and let $\frak q$ be a reduced and irreducible member
of
$|Q|$.
Let $R$ be a line bundle on $X$ such that
\vskip .1 cm
\roster
 \item "3.1.1." $\text H^1(F_i \otimes Q^*)=0$
 \vskip .1 cm
 \item "3.1.2." $\text H ^1(R \otimes \Cal O_\frak q)=0$
 \vskip .1 cm
 \item "3.1.3." $\text H ^1(M_{(F_{i_1} \otimes \Cal O_\frak q)}
\otimes \dots \otimes M_{(F_{i_{q'+1}} \otimes \Cal O_\frak q)}
\otimes R) =0$ for all $0 \leq q' \leq q$
\endroster
Then, for all $-1 \leq q'' \leq q$ and any subset $\{j_k\}
\subseteq \{i\}$ with $\#\{j_k\}=q''+1$ and for all $0 \leq k'
\leq q'' +1$,
$$\text H^1(M_{F_{j_1}} \otimes \dots \otimes M_{F_{j_{k'}}}
\otimes M_{(F_{j_{k'+1}} \otimes \Cal O_\frak q)} \otimes \dots
\otimes M_{(F_{j_{q''+1}}\otimes \Cal O_\frak q)} \otimes R
\otimes \Cal O_\frak q) =0$$
\endproclaim

{\it Proof.}
We prove the result by induction on $q''$. For $q''=-1$ the
corresponding statement is nothing but Condition 3.1.2. Assume that
the result is true for $q''-1$. In order to prove the result
for $q''$ we will use induction on $k'$. If $k'=0$, the
statement is just Condition 3.1.3. Assume that the result is true for
$k' -1$. Because of Condition 3.1.1 we can write for $F_i$ this commutative
diagram:
$$\matrix
&&0&&0&& \\
&& \downarrow && \downarrow && \\
0 &\to& \text H ^0(F_i \otimes Q ^*) \otimes \Cal O _\frak q &
\to & \text H ^0(F_i \otimes Q ^*) \otimes \Cal O _\frak q
&\to& 0 \\
&& \downarrow && \downarrow && \downarrow \\
0& \to & M_{F_i} \otimes \Cal O _\frak q & \to  & \text H ^0(F_i
) \otimes \Cal O_\frak q
&\to&F_i \otimes \Cal O_\frak q & \to & 0 \\
&& \downarrow && \downarrow && \downarrow \\
0& \to& M_{(F_i \otimes \Cal O _\frak q)}&\to&\text H ^0(F_i
\otimes
\Cal O _\frak q) \otimes \Cal O _\frak q& \to &F_i \otimes
\Cal O _\frak q &\to &0 \\
&& \downarrow && \downarrow && \downarrow \\
&&0&&0&&0 \\
\endmatrix$$
Setting $i = j_{k'}$, tensoring the left hand side  vertical
exact sequence by $$M_{F_{j_1}} \otimes \dots \otimes
M_{F_{j_{k'-1}}}
\otimes M_{(F_{j_{k'+1}} \otimes \Cal O_\frak q)} \otimes \dots
\otimes M_{(F_{j_{q''+1}}\otimes \Cal O_\frak q)} \otimes R
\otimes \Cal O_\frak q$$
and taking global sections we obtain this  sequence:
$$\displaylines{\hfill\text H^0(F_{j_{k'}} \otimes Q^*) \otimes
\text H^1(\bigotimes_{r=1}^{k'-1}M_{F_{j_r}} \otimes
\bigotimes_{r=k'+1}^{q''+1} M_{(F_{j_r} \otimes \Cal O_\frak q)}
\otimes R
\otimes \Cal O_\frak q)\hfill\cr
\hfill\to \text
H^1(\bigotimes_{r=1}^{k'}M_{F_{j_r}} \otimes
\bigotimes_{r=k'+1}^{q''+1} M_{(F_{j_r} \otimes \Cal O_\frak q)}
\otimes R
\otimes \Cal O_\frak q)\hfill\cr
\hfill\to \text
H^1(\bigotimes_{r=1}^{k'-1}M_{F_{j_r}} \otimes
\bigotimes_{r=k'}^{q''+1} M_{(F_{j_r} \otimes \Cal O_\frak q)}
\otimes R
\otimes \Cal O_\frak q)\ .\hfill\llap (3.1.4)}$$
The group
$$\text
H^1(M_{F_{j_1}} \otimes \dots \otimes M_{F_{j_{k'-1}}}
\otimes M_{(F_{j_{k'+1}} \otimes \Cal O_\frak q)} \otimes \dots
\otimes M_{(F_{j_{q''+1}}\otimes \Cal O_\frak q)} \otimes R
\otimes \Cal O_\frak q)$$ vanishes by the induction
hypothesis for $q''-1$ and
$$\text
H^1(M_{F_{j_1}} \otimes \dots \otimes M_{F_{j_{k'-1}}}
\otimes M_{(F_{j_{k'}} \otimes \Cal O_\frak q)} \otimes \dots
\otimes M_{(F_{j_{q''+1}}\otimes \Cal O_\frak q)} \otimes R
\otimes \Cal O_\frak q)$$ vanishes by induction on $k'$ (we
have assumed the result to be true for $q''$ and $k'-1$).
Therefore we obtain the vanishing of the group sitting in the middle of
(3.1.4). \boxit{}

\proclaim{Lemma 3.2}
Let $X$ be  a projective variety, let $q$ be a nonnegative integer
and let $F_i$ be a base-point-free line bundle on $X$ for
all $1 \leq i \leq q+1$. Let $Q$ be an effective line
bundle on $X$ and let $\frak q$ be a reduced and irreducible member
of
$|Q|$.
Let $R$ be a line bundle on $X$ such that
\vskip .1 cm
\roster
 \item "3.2.1." $\text H^1(F_i \otimes Q^*)=0$
 \vskip .1 cm
 \item "3.2.2." $\text H ^1(R \otimes Q \otimes \Cal O_\frak
q)=0$
 \vskip .1 cm
 \item "3.2.3." $\text H ^1(M_{(F_{i_1} \otimes \Cal O_\frak q)}
\otimes \dots \otimes M_{(F_{i_{q'+1}} \otimes \Cal O_\frak q)}
\otimes R \otimes Q) =0$ for all $0 \leq q' \leq q$
\endroster
If $\text H^1(M_{F_1} \otimes \dots
\otimes M_{F_{q+1}} \otimes R
) =0 $,
then
$$\text H^1(M_{F_1} \otimes \dots
\otimes M_{F_{q+1}} \otimes R \otimes Q
) =0 \ .$$
\endproclaim

{\it Proof.}
{}From the exact sequence
$$0 \to Q^* \to \Cal O \to \Cal O_\frak q \to 0 \ ,$$
after tensoring by $M_{F_1} \otimes \dots
\otimes M_{F_{q+1}} \otimes R \otimes Q$ and taking global sections
we obtain
$$\displaylines {\text H^1(M_{F_1} \otimes \dots
\otimes M_{F_{q+1}} \otimes R
) \to \text H^1(M_{F_1} \otimes \dots
\otimes M_{F_{q+1}} \otimes R \otimes Q
) \cr
\to \text H^1(M_{F_1} \otimes \dots
\otimes M_{F_{q+1}} \otimes R \otimes Q \otimes \Cal O_\frak q
)\ .}$$
The group $\text H^1(M_{F_1} \otimes \dots
\otimes M_{F_{q+1}} \otimes R
)$ vanishes by hypothesis. To obtain the vanishing of $\text H^1(M_{F_1}
\otimes \dots
\otimes M_{F_{q+1}} \otimes R \otimes Q \otimes \Cal O_\frak q
)$ we use Lemma 3.1 (the line bundle $R$ in Lemma 3.1 is now $R
\otimes Q$ and we set $q'' = q$ and $k' = q+1$). \boxit{}

\proclaim{Lemma 3.3} Let $X$ be a variety. Let $F$, $Q$ and $R$ be line
bundles on $X$ such that $F$ and $F
\otimes Q$ are base-point-free and
$Q$ is effective. Let $\frak q$ be an effective divisor in
$|Q|$, reduced and irreducible. Let $q$ be an integer. Assume that there exists
an integer
$q_0 \leq q$ such that for all $\ 0
\leq l \leq q - q_0 -1$, the following conditions are
satisfied:
\vskip .1 cm
\roster
 \item "3.3.1."  $\text H ^1(F) = \text H ^1(F
\otimes Q ^*)= 0 $
\vskip .1 cm
 \item "3.3.2." $\text H ^1(R
\otimes Q ^{-l} \otimes \Cal O _\frak q)= 0 \ $
\vskip .1 cm
 \item "3.3.3." $\text H ^1(M_{(F
  \otimes \Cal O _\frak q)} ^{\otimes i} \otimes M_{(F
\otimes Q  \otimes \Cal O _\frak q)} ^{\otimes j} \otimes R
\otimes Q ^{-l})= 0 \
\text {for all } \ 1 \leq i+j \leq q-l+1$
\vskip .1 cm
 \item "3.3.4." $\text H ^1(M_{F
} ^{\otimes q-l+1}  \otimes R
\otimes Q ^{-l})= 0$
\vskip .1 cm
\item "3.3.5." $\text H ^1(M_{F
  } ^{\otimes \alpha'} \otimes M_{F
\otimes Q  } ^{\otimes \beta'} \otimes R
\otimes Q ^{-(q-q_0)})= 0 \ \text{for all} \ \alpha' + \beta' q_0 +1 \ .$\endroster

Then, $\text H ^1(M_{F
 } ^{\otimes \alpha} \otimes M_{F \otimes Q} ^{\otimes \beta}
\otimes R \otimes Q^{-m}) = 0$ for all $m$ such that $0 \leq
m \leq q-q_0$  and for all
$\alpha$, $\beta$ nonnegative integers such that $ \alpha +\beta -1 = q-m$. In
particular,
$\text H ^1(M_{F
\otimes Q} ^{\otimes q+1}
\otimes R) = 0$.
\endproclaim

{\it Proof.}
We prove the lemma by induction on $q_0 \leq q'=q-m \leq q$. If
$q-m = q_0$ the conclusion of the theorem is just Condition
3.3.5. Assume that the statement is true for $q' -1 = q-m$. We
will show that it also holds for $q' = q-m$. Now consider $\alpha$
and
$\beta$ such that
$ \alpha +\beta -1
= q'$. We use induction
on $\beta$. If $\beta=0$, the statement is just Condition 3.3.4
considered for $l=q-q'$. Assume that the theorem holds for
$\beta-1$ and we will prove that it holds also for $\beta$. We will
consider two commutative diagrams which yield two exact sequences
relating the bundles
$M_F$, $M_{F \otimes Q}$ and $M_{F \otimes Q \otimes \Cal O
_\frak q}$ (we will set $F \otimes Q \otimes \Cal O
_\frak q = G$, for notational convenience):

$$\matrix
&&0&&0&&0 \\
&& \downarrow && \downarrow &&\downarrow \\
0 &\to& M_F &
\to & \text H ^0(F) \otimes \Cal O_X
&\to& F& \to &0 \\ && \downarrow &&
\downarrow && \downarrow \\ 0& \to &
M_{F \otimes Q} &
\to  & \text H ^0(F
\otimes Q)
\otimes \Cal O_X &\to& F \otimes Q
& \to & 0 \\ && \downarrow &&
\downarrow && \downarrow \\ 0& \to&
K &\to&\text H ^0(G)
\otimes \Cal O_X  & \to &G&\to &0
\\ && \downarrow && \downarrow && \downarrow \\
&&0&&0&&0 \\
\\
\\
&&0&&0 \\
&& \downarrow && \downarrow  \\
0 &\to&\text H^0(G) \otimes Q^*
&
\to & \text H^0(G) \otimes Q^*
&\to&0\\ && \downarrow &&
\downarrow && \downarrow \\ 0& \to &
K & \to  &
\text H^0(G)
\otimes \Cal O_X &\to& G& \to & 0 \\
&& \downarrow && \downarrow && \downarrow \\
0& \to& M_G &\to&\text H^0(G)
\otimes \Cal O_\frak q  & \to &G&\to &0
\\ && \downarrow && \downarrow &&
\downarrow \\ &&0&&0&&0 \\
\endmatrix$$

Note that the exactness at the bottom of the central vertical column of the
first diagram follows from Condition 3.3.1. The two exact sequences we are
interested in are the  ones in the left hand side of each diagram. From the
first one,  after tensoring by
$M_F^{\otimes \alpha}
\otimes M_{F
\otimes Q}
 ^{\otimes \beta-1} \otimes R \otimes Q ^{-(q-q')}$ and taking
global sections, we obtain the sequence
$$\displaylines{\text H^1(M_F^{\otimes \alpha+1} \otimes M_{F \otimes
Q}
 ^{\otimes \beta-1} \otimes R \otimes Q ^{-(q-q')}) \to
\text H^1(M_F^{\otimes \alpha} \otimes M_{F \otimes Q}
^{\otimes \beta} \otimes R \otimes Q ^{-(q-q')}) \cr
\to
\text H^1(K \otimes M_F^{\otimes \alpha} \otimes M_{F \otimes Q}
 ^{\otimes \beta-1} \otimes R \otimes Q ^{-(q-q')})\ .}$$
The group $\text H^1(M_F^{\otimes \alpha+1} \otimes M_{F \otimes
Q}
 ^{\otimes \beta-1} \otimes R \otimes Q ^{-(q-q')})$ vanishes
because, by induction on $\beta$, we have assumed the result to be
true for
$q'=q-m$ and $\beta-1$. Therefore we need only to check that $\text
H^1(K \otimes M_F^{\otimes \alpha} \otimes M_{F \otimes Q}
 ^{\otimes \beta-1} \otimes R \otimes Q ^{-(q-q')})$ vanishes.
For that we use the left hand side exact sequence of the second
diagram. After tensoring it by $M_F^{\otimes \alpha} \otimes M_{F
\otimes Q} ^{\otimes \beta -1} \otimes R \otimes Q^{-(q-q')}$ and
taking global sections we obtain
$$\displaylines{\text H^0(G) \otimes \text H^1(M_F^{\otimes \alpha}
\otimes M_{F
\otimes Q}
^{\otimes \beta-1} \otimes R \otimes Q ^{-(q-q'+1)}) \cr
\to
\text H^1(K \otimes M_F^{\otimes \alpha} \otimes M_{F \otimes Q}
 ^{\otimes \beta-1} \otimes R \otimes Q ^{-(q-q')}) \cr
\to
\text H^1(M_G
\otimes M_F^{\otimes \alpha} \otimes M_{F \otimes Q}
 ^{\otimes \beta-1} \otimes R \otimes Q ^{-(q-q')})}$$
 The group $\text H^1(M_F^{\otimes \alpha}
\otimes M_{F
\otimes Q}
^{\otimes \beta-1} \otimes R \otimes Q ^{-(q-q'+1)})$ vanishes by
induction hypothesis for $q'-1=q-m$. The vanishing of the
cohomology group $\text H^1 (M_G
\otimes M_F^{\otimes \alpha} \otimes M_{F \otimes Q}
 ^{\otimes \beta-1} \otimes R \otimes Q ^{-(q-q')})$ is obtained
from Lemma 3.1 using Conditions 3.3.1 to 3.3.3.

The vanishing of $\text H ^1(M_{F
\otimes Q} ^{q+1}
\otimes R)$ is obtained from the general statement by setting
$m=0$ and $\beta=q+1$. \boxit{}

\vskip .2 cm
 The
remaining lemmas and proposition  are less general
(they are stated for surfaces with geometric genus
$0$) and they  basically yield a slightly more general version of
the vanishing of cohomology obtained in [GP], Proposition 2.1, which we
will need in the arguments of
Sections 4 and 5 when we apply  Lemma 2.8 (concretely when we  check that
Condition 2.8.4 is satisfied).
\proclaim {Lemma 3.4}
Let $X$ be a surface with geometric genus $0$, let $F_1$ and
$F_2$  be two base-point-free, nonspecial line bundles and let $R= F_1
\otimes F_2$. Assume moreover that if $F_1' \equiv F_1$, then
$F_1'$ is base-point-free and nonspecial.
If H$^2(F_2 \otimes {F_1'} ^*) =0$ for all $F_1' \equiv F_1$, then
H$^1(M_R' \otimes {F_1''} ^{\otimes n}) = 0$, for all $n \geq 1 $
 and all $F_1'' \equiv F_1$.
\endproclaim

{\it Proof.}
Mimic  word by word the proof of Lemma 2.5 of [GP] with $R$ playing here the
role of $L$; $F_1$,$F_1'$ and $F_1''$, the role of $B_1$ and
$F_2$, the role of $B_2$.
\boxit{}

\proclaim {Lemma 3.5}
Let $X $ be a surface with geometric genus $0$, let $F_1$ and
$F_2$  be two base-point-free line bundles and let $R= F_1
\otimes F_2$. Assume that $R'$ is nonspecial for all $R' \equiv R$. Assume
also that if $F_1' \equiv F_1$ and $F_2' \equiv F_2$, then $F_1'$ and
$F_2'$ are base-point-free and they satisfy
the conditions   $\text H^1({F_1'}^{\otimes 2})=\text H^1(F_2')=0$
 and $\text H^2(F_2' \otimes {F_1'} ^*) =\text H^2({F_1'}^{\otimes 2}
 \otimes {F_2'} ^*)=0$.

 If $Q$ is any effective  line bundle
 on
$X$ such that either $\text H^1(Q)=0$ or $Q \simeq \Cal
O$,   then
$\text H^1(M_R \otimes R' \otimes Q) = 0$ for any $R' \equiv R$.
\endproclaim

{\it Proof.}
Mimic  word by word the proof of Lemma 2.6 in [GP] with $R$ and $R'$
playing the role of $L$ and  $F_i$ and $F_i'$, the role of
$B_i$).
\boxit{}

\proclaim {Proposition 3.6}
Let $X $ be a surface with geometric genus $0$, let $F_1$ and
$F_2$ be base-point-free line bundles and let $R= F_1
\otimes F_2$. Assume  that if $F_1' \equiv F_1$ and $F_2' \equiv F_2$,
$F_1'$ and $F_2'$ are base-point-free and nonspecial and that they
satisfy
$\text H^2(F_2'
\otimes {F_1'} ^*)=\text H^2(F_1' \otimes {F_2'} ^*)=0$.

Then $\text H ^1(M_R^{\otimes
q+1}\otimes {R'})=0$ for $q = 0,1$  and any $R'
\equiv R$. \endproclaim

{\it Proof.}
Mimic the proof of Proposition 2.1 in [GP], using now Lemmas 3.4 and 3.5
instead of Lemmas 2.5 and 2.6 of [GP] and with $R$ and
$R'$ playing the role of $L$ and  $F_i$ and $F_i'$, the role
of
$B_i$).
\boxit{}

\heading  4. Cohomology vanishings on ruled
elliptic surfaces with invariant $e = -1$  \endheading
In this section $X$ will denote an elliptic ruled
surface with invariant $e=-1$. This means that $X = \bold P (E)$,
where $E$ is a normalized
vector bundle of rank 2  and degree $1$ over a
smooth elliptic curve $C$.
We set $\Cal O(\frak e )= \wedge ^2 E$ and
$e = -\text{deg} \frak e = -1$. We fix a minimal section $C_0$ such that
$\Cal O (C_0) = \Cal O_{\bold P (E)}(1)$. The group Num$(X)$ is
generated by
$C_0$ and by the class of a fiber. We will denote by $f$  the class of a
fiber of $X$. If $\frak a$ is a divisor on $D$, $\frak a f$ will denote the
pullback of $\frak a$ to $X$ by the projection from $X$ to $D$. Sometimes
when deg$\frak a =1$ we will write, by an abuse of notation, $f$ instead of
$\frak a$. The
canonical divisor
$K_X$ is linearly equivalent to $-2C_0 + \frak e f$, and hence
numerically equivalent to $-2C_0 +f$.

As we said at the beginning of Section 2, we want to obtain
sufficient conditions for a line bundle $L$ to satisfy the
property $N_p$. We will be considering $L$ to be a nonspecial,
normally generated line bundle, hence (see Section 1, especially
Lemma 1.2 and (1.2.1)) we are interested in knowing  when
the group
 \text{$\text H ^1(M_L^{\otimes  p' +1}\otimes L)$}
vanishes. Sometimes the approach to a particular problem is
simplified by considering  instead a more general problem. So we  do
here: In  this section and in the next, using the results from
Sections 2 and 3, we obtain sufficient conditions on line bundles
$L_1$ and
$L_2$, so that the group
$\text H ^1(M_{L_1}^{\otimes  p' +1}\otimes L_2)$ vanishes. In
this line, the following are the main results of this section:
\proclaim{Proposition 4.1} Let $B^1 = B_1^1 \otimes \dots \otimes
B_{p+1}^1$ and $B^2 =B_1^2 \otimes \dots \otimes B_{p+1}^2$
be line bundles on $X$ such that $B_i^1 \equiv  B_i^2$ and $B_i^j$ is
in the numerical class of either $2C_0$ or $C_0+f$. Let $P^1$
and $P^2$ be two effective line bundles on $X$ such that $P^1$
is in the numerical class of $aC_0+bf$ for some $a$, $b \geq
0$. If $L_i = B^i \otimes P^i$, then
$$\text H^1 (M_{L_1} ^{\otimes p+1} \otimes L_2)=0$$
\endproclaim

\proclaim{Proposition 4.2}
 Let $p \geq 1$. Let $B^1 B_1^1
\otimes
\dots \otimes
B_{p+1}^1$ and $B^2 = $\linebreak
\text{$ B_1^2\otimes\dots\otimes B_{p+1}^2$}
be line bundles on $X$ such that $B_i^j$ is
in the numerical class of  $2C_0$. Let $P^1$
and $P^2$ be two effective line bundles on $X$. If $L_i = B^i
\otimes P^i$, then
$$\text H^1 (M_{L_1} ^{\otimes p+1} \otimes L_2)=0$$
\endproclaim

To prove the above propositions we
need to do some preliminary work.
We start by recalling when a line bundle on
$X$ is ample, when is base point free, when is effective and when its higher
cohomology vanishes.

\proclaim {Proposition 4.3 ([GP], Proposition 3.1; [Ho1], \S 2,
[Ho2], Proposition 2.3)}

Let $L$ be a line bundle on $X$, numerically equivalent to
$aC_0+bf$. Then

\vskip .1 cm
$$\vbox {\offinterlineskip \hrule
\halign{\vrule
#&\hfil\quad#\quad\hfil&\vrule#&\hfil\quad#\quad\hfil&\vrule#&
\hfil\quad#\quad\hfil&\vrule#&\hfil\quad#\quad\hfil&\vrule#&
\hfil\quad#\quad\hfil&\vrule#
\cr height .1cm &&height .1cm &&height .1cm &&height .1cm&&height
.1cm&&height
.1cm
  \cr
&$a$&&$b$&&$\text h^0(L)$&&$\text h^1(L)$&&$\text h^2(L)$&\cr
height .1cm &&height .1cm &&height .1cm &&height .1cm&&height .1cm
  &&height .1cm\cr
  \noalign{\hrule}
height .1cm &&height .1cm &&height .1cm &&height .1cm&&height
.1cm&&height .1cm
  \cr
  &&&$b > -a/2$&&$>0$&&$0$&&$0$& \cr
height .1cm &&height .1cm &&height .1cm &&height .1cm&&height .1cm
&&height .1cm  \cr
&&&\omit\hrulefill&&\omit\hrulefill&&\omit\hrulefill&
&\omit\hrulefill&\cr
height .1cm &&height .1cm &&height .1cm &&height .1cm&&height .1cm
  &&height .1cm\cr
  &$a \geq 0$&&$b = -a/2$&&$?$&&$?$&&$0$& \cr
height .1cm &&height .1cm &&height .1cm &&height .1cm&&height .1cm
  &&height .1cm \cr
&&&\omit\hrulefill&&\omit\hrulefill&&\omit\hrulefill&
&\omit\hrulefill&\cr
height .1cm &&height .1cm &&height .1cm &&height .1cm&&height .1cm
  &&height .1cm \cr
  &&&$b < -a/2$&&$0$&&$>0$&&$0$& \cr
height .1cm &&height .1cm &&height .1cm &&height .1cm&&height .1cm
  &&height .1cm\cr
\noalign{\hrule}
height .1cm &&height .1cm &&height .1cm &&height .1cm&&height .1cm
  &&height .1cm\cr
  &$a=-1$&&any $b$&&$0$&&$0$&&$0$& \cr
height .1cm &&height .1cm &&height .1cm &&height .1cm&&height .1cm
 &&height .1cm \cr
   \noalign{\hrule}
height .1cm &&height .1cm &&height .1cm &&height .1cm&&height .1cm
  &&height .1cm\cr
  &&&$b > -a/2$&&$0$&&$>0$&&$0$& \cr
height .1cm &&height .1cm &&height .1cm &&height .1cm&&height .1cm
  &&height .1cm\cr
&&&\omit\hrulefill&&\omit\hrulefill&&\omit\hrulefill&
&\omit\hrulefill&\cr
height .1cm &&height .1cm &&height .1cm &&height .1cm&&height .1cm
  &&height .1cm\cr
  &$a \leq -2$&&$b = -a/2$&&$0$&&$?$&&$?$& \cr
height .1cm &&height .1cm &&height .1cm &&height .1cm&&height .1cm
  &&height .1cm\cr
&&&\omit\hrulefill&&\omit\hrulefill&&\omit\hrulefill&
&\omit\hrulefill&\cr
height .1cm &&height .1cm &&height .1cm &&height .1cm&&height .1cm
  &&height .1cm\cr
  &&&$b < -a/2$&&$0$&&$0$&&$>0$& \cr
height .1cm &&height .1cm &&height .1cm &&height .1cm&&height .1cm
  &&height .1cm\cr}  \hrule }$$

\endproclaim

\proclaim {Proposition 4.4 ( [GP], Proposition 3.2)}

$4.4.1.$ There exist three effective line bundles in the numerical
class of $2C_0 -f$. They are $\Cal O (2C_0 + (\frak e + \eta
_i)f)$, where the $\eta _i$s are the nontrivial degree $0$
divisors corresponding to the three nonzero torsion points
in Pic$^0(C)$. The unique element in $|2C_0 + (\frak e + \eta
_i)f)|$ is  a smooth elliptic curve $Y_i$

$4.4.2.$ For each $n > 0$, there are only four effective line
bundles numerically equivalent to $n(2C_0 -f)$. They are $\Cal O
(2nC_0 + n(\frak e + \eta
_i)f)$ and $\Cal O (2nC_0 + n\frak e f)$. The only smooth
(elliptic) curves (and indeed the only irreducible curves) in these
numerical classes are general members in $|4C_0 + 2\frak e f|$.

The number of linearly independent global sections of these line
bundles are summarized in the following table:
\vskip .2 cm
$$\vbox {\offinterlineskip \hrule
\halign{\vrule
#&\hfil\quad#\quad\hfil&\vrule#&\hfil\quad#\quad\hfil&\vrule#&
\hfil\quad#\quad\hfil&\vrule#&\hfil\quad#\quad\hfil&\vrule#&
\hfil\quad#\quad\hfil&\vrule#&
\hfil\ #\ \hfil&\vrule#&
\hfil\quad#\quad\hfil&\vrule#&
\hfil\quad#\quad\hfil&\vrule#
\cr
height .1cm &&height .1cm &&height .1cm &&height .1cm&&height
.1cm&&height .1cm&&height
.1cm&&height .1cm&&height .1cm
  \cr
  &$n \geq 0$&&$0$&&$1$&&$2$&&$3$&&$\dots$&&$2m$&&$2m+1$& \cr
height .1cm &&height .1cm &&height .1cm &&height .1cm&&height .1cm
&&height .1cm &&height .1cm&&height .1cm&&height .1cm \cr
\noalign{\hrule}
height .1cm &&height .1cm &&height .1cm &&height .1cm&&height .1cm
  &&height .1cm&&height .1cm&&height .1cm&&height .1cm\cr
  &$\text h ^0(\Cal O (2nC_0 + n \frak e
f))$&&$1$&&$0$&&$2$&&$1$&&&&$m+1$&&$m$& \cr
height .1cm &&height .1cm &&height .1cm &&height .1cm&&height .1cm
  &&height .1cm &&height .1cm&&height .1cm&&height .1cm\cr
\noalign{\hrule}
height .1cm &&height .1cm &&height .1cm &&height .1cm&&height .1cm
  &&height .1cm &&height .1cm&&height .1cm&&height .1cm\cr
  &$\text h ^0(\Cal O (2nC_0 + n (\frak e + \eta _i)
f))$&&$0$&&$1$&&$1$&&$2$&&&&$m$&&$m+1$& \cr
height .1cm &&height .1cm &&height .1cm &&height .1cm&&height .1cm
  &&height .1cm&&height .1cm&&height .1cm&&height .1cm\cr }
\hrule }$$ \endproclaim
\vskip .2 cm
(4.4.3) We will fix one of the smooth elliptic curves in the numerical class
of $2C_0-f$ and we will call it $E$.
\vskip .2 cm
\proclaim {Proposition 4.5 ([H], V.2.21.b; [GP], Proposition 3.5 and
Remark 3.5.3)}\linebreak
Let $L$ be a line bundle on $X$ in the numerical class of $aC_0 +
bf$.
\roster
\item"4.{5}.1" $L$ is ample iff $a >0$ and $b >-\frac 1 2 a$
\item"4.{5}.2" $L$ is base point free if $a \geq 0$,
$a+b \geq 2$ and 	$a+2b \geq 2$.
\item"4.{5}.3" $L$ is ample and base point free iff $a \geq 1$,
$a+b \geq 2$ and 	$a+2b \geq 2$.
\endroster
\endproclaim

We will need some lemmas dealing with the vanishing of the cohomology of
certain bundles on curves:

\proclaim{Lemma 4.6}
Let $p \geq -1$ and
let $B_i$ be a line bundle on $\bold P^1$ for all
$1 \leq i \leq p+1$, such that $b_i \text{deg}(B_i) \geq 1$. Let $L$ be a line bundle on $\bold P^1$
such that $l = \text{deg}(L) \geq p$. Then $\text H
^1(M_{B_1} \otimes \dots \otimes M_{B_{p+1}}
\otimes L) =0$
\endproclaim

{\it Proof.}
Note in the first place that each $B_i$ is base-point-free so
it makes sense to define $M_{B_i}$. The bundle $M_{B_i}$ is
isomorphic to $\Cal O_{\bold P^1}(-1)^{\oplus b_i}$ (the
sequence defining $M_{B_i}$ is the sheafification of
$$0 \to S^{b_i}(-1) \to S^{b_i+1} \to S(b_i) \to 0$$
where $S$ denotes the homogeneous coordinate ring of $\bold
P^1$). Hence using that $$\text H
^1(M_{B_1} \otimes \dots \otimes M_{B_{p+1}}
\otimes L) = \text H^1(\Cal O(l-p-1)^{\oplus (\prod _{i=1}^{p+1}b_i)})$$
and that $l-p-1 \geq -1$, we obtain the result. \boxit{}

\proclaim{Lemma 4.7} Let $Y$ be a smooth elliptic
curve.
Let $p \geq -1$ and let $B_i$ be a  line bundle on $Y$
for all
$1 \leq i \leq p+1$. Let $L$ be another line bundle on $Y$.
Let $b_i = \text{deg}(B_i) \geq 2$ and let $l = \text{deg}(L)$. If  $\sum
_{i=1}^{p+1} \frac {b_i} {b_i -1} < l$,
then
$$\text H^1(M_{B_1} \otimes \dots \otimes M_{B_{p+1}} \otimes L) =0
\ .$$ In particular, if $b_i \geq p+3$ for all $1 \leq i \leq p+1$
and $l \geq p+2$, then $$\text H^1(M_{B_1} \otimes \dots
\otimes M_{B_{p+1}}
\otimes L) =0 \ .$$
\endproclaim

{\it Proof.}
Note first that the $B_i$s are base-point-free, since their
degrees are greater or equal than $2$; hence $M_{B_i}$ makes sense.
Let
$r_i = r(B_i) = h^0(B_i) - 1$. Then $$\displaylines{\text{rk}(M_{B_1}
\otimes \dots
\otimes M_{B_{p+1}}
\otimes L) = r_1
\cdots r_{p+1} \quad  \text{and} \cr
\text{deg}(M_{B_1} \otimes \dots \otimes M_{B_{p+1}}
\otimes L) = l \cdot r_1 \cdots r_{p+1} - \sum _{i=1} ^{p+1} b_i
\cdot  r_1
\cdots
\hat{r_i} \cdots r_{p+1}  \cr}$$
and therefore $$\mu (M_{B_1}
\otimes \dots \otimes M_{B_{p+1}} \otimes L)
=l - \sum _{i=1} ^{p+1} \frac {b_i}
{r_i}$$
The bundle $M_{B_1}
\otimes \dots \otimes M_{B_{p+1}} \otimes L$ is semistable by [Bu],
Theorem
1.2 and [Mi], Corollary {4.9} and \S 5. Therefore $\text
H^1(M_{B_1}
\otimes \dots \otimes M_{B_{p+1}} \otimes L) =0$ if $\mu(M_{B_1}
\otimes \dots \otimes M_{B_{p+1}} \otimes L) >2g(Y) -2$.
Since
$Y$ is elliptic, $r_i=b_i-1$ and $2g(Y)-2 =0$ and the conclusion of
the lemma is clear.
\boxit{}
\vskip .2 cm
Now we prove some lemmas which have in account the particular
properties of elliptic ruled surfaces.
\proclaim{Lemma 4.8}
Let $L =B_1 \otimes B_2$
and $L' = B_3$ be line bundles on $X$ satisfying the following
properties:
\roster
\item " 4.{8}.1" $B_1 \equiv B_3$
\item " 4.{8}.2" $B_1$ is in the numerical class of
$2C_0$ or in the numerical class of $C_0 +f$.
\item "4.{8}.3" $B_2$ is the numerical class of $2C_0$, in the
numerical class of $C_0 +f$ or in the numerical class of $2f$
\endroster
Let $P$ and $P'$ be effective. Then the map
$$\text H ^0 (L \otimes P) \otimes \text H^0(L' \otimes P')
@>\alpha >> \text H^0(L \otimes L' \otimes P \otimes P')$$
surjects.
\endproclaim

{\it Proof.}  We start by noting that the multiplication
map $$\text H ^0(L \otimes P) \otimes \text H^0(L') @>\beta >>
\text H^0(L
\otimes L' \otimes P)$$ surjects. This is a consequence of Theorem
1.3. Indeed.  The line bundle $L'$ is base-point-free by
Proposition {4.5} and since
$L \otimes P
\otimes {L'}^*
\equiv B_2 \otimes P$ and $L \otimes P \otimes {L'}^{-2} \equiv
B_2
\otimes B_3 ^* \otimes P$,
it follows from Proposition {4.3} that $\text H^1(L \otimes P
\otimes {L'}^*
) = \text H^2(L \otimes P \otimes {L'}^{-2}) =0$. Since, by
Proposition {4.3}, $\text H^1(L')=0$, the surjectivity of $\beta$ is
equivalent to the vanishing of the group $\text H^1(M_{L \otimes P}
\otimes L')$. Since, again by
Proposition {4.3}, $\text H^1(L' \otimes P')=0$,  the surjectivity of
$\alpha$ is equivalent to the vanishing of $\text H^1(M_{L \otimes
P} \otimes L' \otimes P')$.  We
use Lemma 3.2 to prove the latter vanishing. We can assume without
loss of generality that $P' \Cal O(aC_0 + bf + cE)$.
Thus we carry out induction on $(a,b,c)$. If $(a,b,c) (0,0,0)$ the content of the statement we want to prove is
nothing but the vanishing of $\text H^1(M_{L \otimes P}
\otimes L')$, which we have just shown. Now we assume the
result to be true for $(a-1,0,0)$ and we will prove that it is also
true for $(a,0,0)$. For that we apply Lemma 3.2 to $q=0$,
$F_1=L \otimes P$, $Q=\Cal O(C_0)$,
$\frak q=C_0$ and $R=L' \otimes \Cal
O((a-1)C_0)$.  We need to see that the conditions
required by  Lemma 3.2 are satisfied. For Condition
3.2.1 it is enough to check that
$\text H^1(L \otimes P \otimes \Cal O(-C_0))=0$ and this is true
by Proposition {4.3}. Using that deg$(L' \otimes \Cal
O_{C_0}(aC_0)) \geq 3 >0$ we see that Condition 3.2.2 is satisfied. Condition
3.2.3 follows from Lemma 4.7 because deg$(L' \otimes \Cal
O_{C_0}(aC_0)) \geq 3 $ and deg$(L \otimes P \otimes \Cal O_{C_0})
\geq 4$. The argument for the induction on $b$ and $c$
is analogous (in the case of $b$ we use Lemma 4.6
instead of Lemma 4.7) to the one we have just made and we will not
show it here. \boxit{}

\proclaim{Lemma 4.9}   Let $a$, $b$ be two integers such that $a
\geq 1$,
$a+b\geq 4$ and $a+2b \geq 4$. Let
$L$ be a line bundle in the numerical class of
$aC_0+bf$ and $P$ a line bundle whose numerical
class contains an effective representative. Then
$$\text H^2(M_L^{\otimes 2} \otimes P) =0$$ \endproclaim

{\it Proof.} Note first that $L$ is base-point-free
(c.f. Proposition {4.5}) and therefore it make sense to talk
about $M_L$. From exact sequence 1.1 we obtain these two
exact sequences:
$$ \displaylines {\hfill\text H^1(M_L \otimes L \otimes P)
\to \text H^2(M_L^{\otimes 2} \otimes P) \to \text
H^0(L) \otimes \text H^2(M_L \otimes P)\hfill\llap{\text{and}\quad} \cr
\hfill\text H^1(L \otimes P) \to \text H^2(M_L \otimes P) \to
\text H^0(L) \otimes \text H^2(P)\hfill \cr}$$
The vanishing of $\text H^1(M_L \otimes L \otimes P)$
follows from Lemma  4.5 and Proposition 4.3. The vanishings of
$\text H^1(L
\otimes P)$ and $\text H^2(P)$ follow from Proposition {4.3}, and
hence we obtain the result. \boxit{}
\vskip .1 cm

\proclaim {Lemma 4.10} Let $B$ be a line bundle in the numerical class of
$2(p+1)C_0$ for some $p \geq 1$. Then one can choose a divisor $\frak d$ of
degree $1$ on
$C$ and $B_i$s in the numerical class of $2C_0$ for all $1 \leq i
\leq p+1$ such that $B_i \otimes \Cal O(-\frak d f)$ is effective
for all $1 \leq i
\leq p$, but $B
\otimes
\Cal O(-(p+1)\frak d f)$ is not effective.
\endproclaim

{\it Proof.}
The line bundle $B$ is equal to $\Cal O(2(p+1)C_0 + \frak a f)$ for
some
degree $0$ divisor $\frak a$ on $C$. Choose $\frak d$ satisfying
$2(\frak a +(p+1)(\frak e - \frak d)) \not\sim 0$ and set, for all $1 \leq i
\leq p$,
$B_i$ equal to
$\Cal O(2C_0 + (\frak d -\frak e - \eta)$ for some divisor $\eta$
on
$C$ such that
$2\eta
\sim 0$ and $\eta \not\sim 0$. Then Proposition {4.4} implies that
$B_i \otimes
\Cal O(-\frak d f)$ is effective for all $1 \leq i
\leq p$ and that $B
\otimes
\Cal O(-(p+1)\frak d f)$ is not effective. \boxit{}
\vskip .2 cm
We are now ready to give the proof of Propositions 4.3 and 4.4:

(4.11) {\it Proof of Proposition 4.1.}

{\it Step 1.} $\text H^1 (M_{B^1} ^{\otimes p+1} \otimes L_2)=0$

We will use Lemma 2.8.  The set
$\frak B$ will consist of those line bundles on $X$ belonging
either to the numerical class of $2C_0$ or to the numerical class
of
$C_0+f$; the set $\frak P$ will be the set of all effective line
bundles of
$X$ and
$q_0$ will be equal to
$1$. Therefore if
$\frak B$ and $\frak P$ satisfy the conditions of Lemma 2.8, we are
done (simply take
$n$ to be $0$ and $Q = P^2$  in the conclusion of
Lemma 2.8). Conditions 2.8.1 and 2.8.2 are satisfied, (c.f.
Proposition {4.3} and Proposition {4.5}). For Condition 2.8.3 note
that
$B_1
\otimes B_2
\otimes B_3^*$  is an effective line bundle for any $B_1$,
$B_2$ and
$B_3$ in $\frak B$ (this follows again from Proposition {4.3}). On
the other hand, the line bundle $R_3$ defined in the statement
of Lemma 2.8 satisfies the hypothesis for the line bundle $L$ in
Lemma {4.9}. Thus applying the mentioned lemma we are done. For
Condition 2.8.4 we have to show that
$$\text H^1(M_{B_1 \otimes B_2 \otimes C_1 \otimes
\dots \otimes C_n}^{\otimes 2} \otimes B_1' \otimes B_2'
\otimes P) = 0$$
for all $B_1, B_2, B_1', B_2', C_1, \dots ,C_n
\in
\frak B$ and $P$ effective line bundle satisfying the condition $B_i \equiv
B_i'$. Note that $C_1
\otimes \dots \otimes C_n$ is numerically equivalent to $\Cal O(aC_0 +bf)$ for
some
$a, b \geq 0$, so we will prove  this more general fact instead:
 Let $P^1$ be an effective line bundle in the numerical class of
$a_1C_0+b_1f$ for some $a_1,b_1 \geq 0$ and let $P^2$ be another
effective line bundle. Then
$$\text H^1(M_{B_1^1 \otimes B_2^1 \otimes P^1} ^{\otimes 2}
\otimes B_1'
\otimes B_2' \otimes P^2) =0 \ .\leqno (4.11.1)$$
First we show
using Lemma 3.1 that
$\text H^1(M_{B_1 \otimes B_2 \otimes P^1}^{\otimes 2} \otimes B_1' \otimes
B_2' ) = 0$. By Lemma {4.7} we may assume without loss of generality
that
\vskip .1 cm
({4.11}.2)
 $P^1 = \Cal O(a_1C_0 +b_1f)$, and if $B_1' \equiv B_2'
\equiv \Cal O(2C_0)$, the line bundle $B_1' \otimes \Cal O(-f)$ is
effective and the line bundle
$B_1' \otimes B_2'
\otimes \Cal O(-2f)$ is  not effective. In particular,
$\text H^1(B_1' \otimes B_2'
\otimes \Cal O(-2f))=0$.
\vskip .1 cm
We use induction on $(a_1,b_1)$. If
$(a_1,b_1) = (0,0)$, the result follows from Proposition 3.6.
Now assume that
$$\text
H^1(M_{B_1 \otimes B_2 \otimes \Cal O((a_1-1)C_0)}^{\otimes 2}
\otimes B_1' \otimes B_2') =0$$
for $a_1 \geq 1$. We apply Lemma 3.3 to $F = B_1 \otimes B_2 \otimes
\Cal O((a_1-1)C_0)$, $Q=\Cal O(C_0)$, $\frak q = C_0$, $R=B_1'
\otimes B_2'$, $q=1$, $q_0=-1$, $\alpha=0$ and $m=0$.
Condition 3.3.1 is satisfied by Proposition {4.3}. Condition 3.3.2
is
satisfied because deg$(B_1' \otimes B_2' \otimes \Cal
O_{C_0}(-lC_0)) \geq 3$ for $l = 0,1$. We check that Condition 3.3.3 is
satisfied by using Lemma 4.7, noting that $$\text{deg}(B_1 \otimes
B_2
\otimes \Cal
O_{C_0}(a_1C_0)) > \text{deg}(B_1 \otimes B_2 \otimes \Cal
O_{C_0}((a_1-1)C_0)) \geq 4 $$
and that
$$\text{deg}(B_1' \otimes
B_2' \otimes \Cal O_{C_0}(-lC_0)) \geq 3 \ .$$
Condition 3.3.4 requires that $$\text H^1(M_{B_1 \otimes B_2
\otimes \Cal O((a_1-1)C_0)} \otimes B_1' \otimes B_2' \otimes
\Cal O(-C_0))=0 \ ,$$ which is a consequence of Lemma 4.8 and
Proposition 4.3, and that
$$\text H^1(M_{B_1 \otimes B_2
\otimes \Cal O((a_1-1)C_0)}^{\otimes 2} \otimes B_1' \otimes B_2'
)=0$$
which is true by the induction hypothesis on $a_1-1$.
Condition 3.3.5 requires the vanishing of $\text H^1(B_1'
\otimes B_2' \otimes \Cal O(-2C_0)) =0$ which follows from
Proposition {4.3}.
Now  we carry out induction on $b_1$. If
$b_1 =0$, the required statement has just been proven. Assume that
the result is true for $b_1 -1$  ($b_1 \geq 1$). We will use
again Lemma 3.3 setting $F = B_1 \otimes B_2
\otimes \Cal O(a_1C_0+(b_1-1)f)$, \text{$R=B_1'
\otimes B_2' $,} $Q = \Cal O(f)$, $\frak q = f$,
$q=1$ and $q_0=-1$. Condition 3.3.1 is satisfied because
of Proposition {4.3}. Condition 3.3.2 is satisfied because deg$(B_1'
\otimes B_2'
\otimes \Cal O_f(-lf)) \geq 2$. Condition 3.3.3
follows from Lemma 4.6, since
$$\displaylines{\hfill\text{deg}(B_1
\otimes B_2
\otimes \Cal O_f(a_1C_0+b_1f)) \geq 2 \hfill \cr
\hfill\text{deg}(B_1 \otimes B_2
\otimes \Cal O_f(a_1C_0+(b_1-1)f)) \geq 2 \hfill\llap{\text{and} \quad}\cr
\hfill\text{deg}(B_1' \otimes B_2'
\otimes \Cal O_f(-lf)) \geq 2 \hfill\cr}$$
Condition 3.3.4 follows by induction hypothesis on $b_1-1$ and
from Lemma 4.8 and ({4.11}.2). Condition 3.3.5 follows
from Proposition {4.3} and ({4.11}.2).

To finish the proof of ({4.11}.1) we apply Lemma 3.2
inductively (as done for instance in the proof of Lemma 4.8) setting
the line bundles
$F_1$ and
$F_2$  both equal to
$B_1
\otimes B_2
\otimes P^1$, $Q$ equal
to $\Cal O(C_0)$, $\Cal O(f)$ or $\Cal O(E)$, $\frak q$
equal to $C_0$, $f$, or $E$ and $q = 2$.
\vskip .3 cm
{\it Step 2.} $\text H^1(M_{L_1} ^{\otimes p+1} \otimes L_2)
=0$.

Again by Lemma {4.10} we may assume without loss of generality the
following:
\vskip .2 cm
({4.11}.3) $L_1 = B_1^1
\otimes \dots \otimes B_{p+1}^1 \otimes \Cal O(a_1C_0 + b_1f)$
and if all the $B_i^2$s are in the
numerical class of $2C_0$, then $B_i^2 \otimes
\Cal O(-f)$ is effective for all $2 \leq i \leq p+1$ but
$B_1^2 \otimes \dots \otimes B_{p+1}^2 \otimes
\Cal O(-pf)$ is not effective, for all $p \geq 1$. In particular
$\text H^1(B_1^2 \otimes \dots \otimes B_{p+1}^2 \otimes
\Cal O(-pf)) =0$.
\vskip .2 cm
We will prove that
$$  \text H^1(M_{B_1^1 \otimes \dots \otimes B_{p'+1}^1 \otimes P^1} ^{\otimes
p'+1}
\otimes B_1^2
\otimes \dots \otimes B_{p'+1}^2 \otimes P^2) =0 \leqno ({4.11}.4)$$
for all $1 \leq p' \leq p$. We use induction on $p'$.
If
$p'=1$ we must prove that
$$\text H^1(M_{B_1^1 \otimes B_2^1 \otimes P^1} ^{\otimes 2}
\otimes B_1^2
\otimes B_2^2 \otimes P^2) =0 \ .$$
This is the content of ({4.11}.1).

Now we assume that ({4.11}.4) holds for $1, \dots , p'-1$ ($p'
\geq 2$) and we prove that it holds also for $p'$. Again we make
induction on $(a_1, b_1)$. If $(a_1,b_1) = (0,0)$ the
statement was proven in Step 1. Assume that the result is true for
$(a_1-1,0)$. We apply Lemma 3.3 to $F=B^1 \otimes \Cal
O((a_1-1)C_0)$, $R = L_2$, $\frak q = C_0$, $q=p'$ and $q_0=-1$.
Condition 3.3.1 is satisfied by Proposition {4.3}. Condition 3.3.2
follows from the fact that deg$(L_2 \otimes \Cal
O_{C_0}(-lC_0)) \geq p'+2 \geq 4$ for all $0 \leq l \leq p'$. Condition 3.3.3
follows from Lemma 4.7 and from the fact that  deg$(L_2 \otimes \Cal
O_{C_0}(-lC_0)) \geq p'+2$ and $$\text{deg}(B^1 \otimes \Cal
O_{C_0}(a_1C_0)) \geq \text{deg}(B^1 \otimes \Cal
O_{C_0}((a_1-1)C_0)) \geq 2p'+2 \geq p'+3 \ .$$
Condition 3.3.4 requires the vanishing of
$$\text H^1(M_{B^1
\otimes \Cal
O((a_1-1)C_0)} ^{\otimes p''+1} \otimes L_2 \otimes \Cal O(-lC_0))
\eqno ({4.11}.5)$$ for
$l=p'-p''$ and $0 \leq p'' \leq p'$. If $p'' = p'$, the
vanishing of ({4.8}.5) is simply the induction hypothesis for
$a_1-1$. If $1 \leq p''
\leq p'-1$, the vanishing of ({4.11}.5) follows from the induction
hypothesis on
$1, \dots , p'-1$. Indeed. The line bundle  $B^1 \otimes \Cal
O((a_1-1)C_0)$ can be written as the tensor product of $B_1^1
\otimes
\dots
\otimes B_{p''+1}^1$ with an effective line bundle numerically equivalent to
$aC_0+bf$ for some $a,b \geq 0$. The line bundle $L_2 \otimes \Cal
O(-lC_0)$ can be written as the tensor product of $B_1^2 \otimes
\dots \otimes B_{p''+1}^2$ with an effective line bundle, since
$B_i^2 \otimes \Cal O(-C_0)$ is effective. If $p'' =0$, the
vanishing of ({4.11}.5) follows from Lemma 4.8 and Proposition 4.3.
Condition 3.3.5 requires the vanishing of
$\text H^1(L_2 \otimes \Cal O(-(p'+1)C_0))$ which follows from
Proposition {4.3}.

The induction argument on $b_1$ is similar to the one on $a_1$
and we will only highlight here the differences and the delicate
points. We make again iterated use of Lemma 3.3. Condition
3.3.3 follows from Lemma 4.3. Condition 3.3.4 is
obtained as before (assumption  ({4.11}.3) assures us that
$B_i^2 \otimes \Cal O(-f)$ is effective for all $1 \leq i \leq p'$).
 Condition
3.3.5 is obtained from Propositions {4.3} and {4.4} and assumption
({4.11}.3).
\boxit{}

(4.12)
{\it Proof of Proposition 4.2.} Without loss of generality we may
assume that
$P^1$ is isomorphic
to
$\Cal O(a_1C_0 +b_1f +c_1E)$. We prove the result by induction
on $p$. First we prove it for $p=1$. We will use induction on
$c_1$. If
$c_1=0$ the result follows from Proposition {4.1}. Assume that the
result is true for
$(a_1,b_1,c_1-1)$ and $c_1 \geq 1$. We apply Lemma 3.3
to \text{$F B^1
\otimes \Cal O(a_1C_0 +b_1f +(c_1-1)E)$,} $\frak q = E$,
$R = L_2$, $q =1$ and $q_0=-1$. We have to check that the
conditions of Lemma 3.3 are satisfied.  Condition 3.3.1 follows
from Proposition {4.3}. Condition 3.3.2 follows from the fact that
$\text{deg}(L_2 \otimes \Cal O_E(-lE)) = \text{deg}(L_2 \otimes
\Cal O_E)>0$. Condition 3.3.3 follows from Lemma 4.7 using
the fact that \text{$\text{deg}(L_2 \otimes \Cal O_E(-lE))
\geq 4 $} and $\text{deg}(B^1
\otimes \Cal O(a_1C_0 +b_1f +(c_1-1)E)) \geq 4$. Condition
3.3.4 requires the vanishing of $\text H^1(M_{B^1
\otimes \Cal O(a_1C_0 +b_1f +(c_1-1)E)} \otimes L_2 \otimes
\Cal O(-E))$ which follows from Lemma 4.8 and the
vanishing of $\text H^1(M_{B^1
\otimes \Cal O(a_1C_0 +b_1f +(c_1-1)E)}^{\otimes 2} \otimes L_2
)$ which follows from the induction hypothesis on $c_1-1$.
Condition 3.3.5 follows from Proposition {4.3}.

Now let us assume the result to be true for $1, \dots,
p-1$. To prove the result for $p \geq 2$ we will again use
induction on
$c_1$. If
$c_1=0$ the result follows from Proposition {4.11}. Assume that the
result is true for
$(a_1,b_1,c_1-1)$ and $c_1 \geq 1$. We apply Lemma 3.3 to $F B^1
\otimes \Cal O(a_1C_0 +b_1f +(c_1-1)E)$, $\frak q = E$,
$R = L_2$, $q=p$ and $q_0=-1$. We see now that the conditions
of Lemma 3.3 are satisfied. Condition 3.3.1 follows
from Proposition {4.3}. Condition 3.3.2 follows from the fact that
$\text{deg}(L_2 \otimes \Cal O_E(-lE)) = \text{deg}(L_2 \otimes
\Cal O_E)>0$. Condition 3.3.3 follows from Lemma 4.7 using
the fact that $$\displaylines{\text{deg}(L_2 \otimes \Cal O_E(-lE))
\geq 2p+2 > p+2 \qquad  \text{and} \cr
\text{deg}(B^1
\otimes \Cal O_E(a_1C_0 +b_1f +(c_1-1)E)) \geq 2p+2 \geq
p+3 \ .\cr}$$
Condition 3.3.4 requires the vanishing of $$\text H^1(M_{B^1
\otimes \Cal O(a_1C_0 +b_1f +(c_1-1)E)} ^{\otimes p'+1} \otimes
L_2
\otimes
\Cal O(-lE))$$ for all $0 \leq p' \leq p$ and $l =p
-p'$. If $p'=p$ the vanishing follows from the
induction hypothesis on $c_1-1$. If $1 \leq p' \leq
p-1$ the vanishing follows from the induction
hypothesis on $1, \dots, p-1$. Finally, if $p'=0$ the
vanishing follows from Lemma 4.8 and Proposition 4.3. Condition
3.3.5 requires the vanishing of $\text H^1(L_2 \otimes
\Cal O(-(p+1)E))$ which follows from Proposition {4.3}. \boxit{}

\heading 5. Cohomology vanishings on
elliptic ruled surfaces with invariant $e \geq
0$  \endheading
In this section we duplicate for an elliptic ruled surface of invariant $e
\geq 0$ the work done in the previous one for an elliptic ruled surface with
invariant $e=-1$.
Thus $X$ will denote throughout this section an elliptic ruled
surface with invariant $e \geq 0$. Again $C_0$
will be a minimal section of $X$. We will denote by $f$ the class of a
fiber of $X$. If $\frak a$ is a divisor on $D$, $\frak a f$ will denote the
pullback of $\frak a$ to $X$ by the projection from $X$ to $D$. Sometimes
if $\text{deg}\ \frak a =1$ we will write, by an abuse of notation, $f$
instead of
$\frak a f$. The canonical divisor
$K_X$ is linearly equivalent to $-2C_0 + \frak e f$, and hence
numerically equivalent to $-2C_0 -ef$.

Our main result in this section is
\proclaim{Proposition 5.1} Let $B^1 = B_1^1 \otimes \dots
\otimes B_{p+1}^1$ and $B^2 =B_1^2 \otimes \dots \otimes
B_{p+1}^2$ be line bundles such that
$B_i^j$ is in the numerical class of $C_0+(e+2)f$.
Let $P^1$ and $P^2$ be two effective line bundles on $X$ such
that $P^j$ is in the numerical class of $a_j(C_0+ef)+b_jf$ for
some
$a_j$,
$b_j \geq 0$. If $L_i = B^i \otimes P^i$, then
$$\text H^1 (M_{L_1} ^{\otimes p+1} \otimes L_2)=0$$
\endproclaim

Before we prove Proposition 5.1 we need to recall some properties
of the line bundles on $X$.

\proclaim {Proposition 5.2 ( [GP] Proposition 3.1 or [Ho1], \S 2)}

Let $L$ be a line bundle on $X$, numerically equivalent to
$aC_0+bf$. Then

\vskip .1 cm
$$\vbox {\offinterlineskip \hrule
\halign{\vrule
#&\hfil\quad#\quad\hfil&\vrule#&\hfil\quad#\quad\hfil&\vrule#&
\hfil\quad#\quad\hfil&\vrule#&\hfil\quad#\quad\hfil&\vrule# \cr
height .1cm &&height .1cm &&height .1cm &&height .1cm&&height .1cm
  \cr
&$a$&&$b$&&$\text h^0(L)$&&$\text h^2(L)$&\cr
height .1cm &&height .1cm &&height .1cm &&height .1cm&&height .1cm
  \cr
  \noalign{\hrule}
height .1cm &&height .1cm &&height .1cm &&height .1cm&&height .1cm
  \cr
  &&&$b > 0$&&$>0$&&$0$& \cr
height .1cm &&height .1cm &&height .1cm &&height .1cm&&height .1cm
  \cr    &&&\omit\hrulefill&&\omit\hrulefill&&\omit\hrulefill&\cr
height .1cm &&height .1cm &&height .1cm &&height .1cm&&height .1cm
  \cr
  &$a \geq 0$&&$b = 0$&&$?$&&$0$& \cr
height .1cm &&height .1cm &&height .1cm &&height .1cm&&height .1cm
  \cr &&&\omit\hrulefill&&\omit\hrulefill&&\omit\hrulefill&\cr
height .1cm &&height .1cm &&height .1cm &&height .1cm&&height .1cm
  \cr
  &&&$b < 0$&&$0$&&$0$& \cr
height .1cm &&height .1cm &&height .1cm &&height .1cm&&height .1cm
  \cr
\noalign{\hrule}
height .1cm &&height .1cm &&height .1cm &&height .1cm&&height .1cm
  \cr
  &$a=-1$&&any $b$&&$0$&&$0$& \cr
height .1cm &&height .1cm &&height .1cm &&height .1cm&&height .1cm
  \cr
   \noalign{\hrule}
height .1cm &&height .1cm &&height .1cm &&height .1cm&&height .1cm
  \cr
  &&&$b > -e$&&$0$&&$0$& \cr
height .1cm &&height .1cm &&height .1cm &&height .1cm&&height .1cm
  \cr    &&&\omit\hrulefill&&\omit\hrulefill&&\omit\hrulefill&\cr
height .1cm &&height .1cm &&height .1cm &&height .1cm&&height .1cm
  \cr
  &$a \leq -2$&&$b = -e$&&$0$&&$?$& \cr
height .1cm &&height .1cm &&height .1cm &&height .1cm&&height .1cm
  \cr &&&\omit\hrulefill&&\omit\hrulefill&&\omit\hrulefill&\cr
height .1cm &&height .1cm &&height .1cm &&height .1cm&&height .1cm
  \cr
  &&&$b < -e$&&$0$&&$>0$& \cr
height .1cm &&height .1cm &&height .1cm &&height .1cm&&height .1cm
  \cr}  \hrule }$$

$$\vbox {\offinterlineskip \hrule
\halign{\vrule
#&\hfil\quad#\quad\hfil&\vrule#&\hfil\quad#\quad\hfil&\vrule#&
\hfil\quad#\quad\hfil&\vrule# \cr
height .1cm &&height .1cm &&height .1cm &&height .1cm
  \cr
&$a$&&$b$&&$\text h^1(L)$&\cr
height .1cm &&height .1cm &&height .1cm &&height .1cm
  \cr
  \noalign{\hrule}
height .1cm &&height .1cm &&height .1cm &&height .1cm
  \cr
  &&&$b > ae$&&$0$& \cr
height .1cm &&height .1cm &&height .1cm &&height .1cm
  \cr    &&&\omit\hrulefill&&\omit\hrulefill&\cr
height .1cm &&height .1cm &&height .1cm &&height .1cm
  \cr
  &$a \geq 0$&&$b = ae$&&$?$& \cr
height .1cm &&height .1cm &&height .1cm &&height .1cm
  \cr &&&\omit\hrulefill&&\omit\hrulefill&\cr
height .1cm &&height .1cm &&height .1cm &&height .1cm
  \cr
  &&&$b < ae$&&$>0$& \cr
height .1cm &&height .1cm &&height .1cm &&height .1cm
  \cr
\noalign{\hrule}
height .1cm &&height .1cm &&height .1cm &&height .1cm
  \cr
  &$a=-1$&&any $b$&&$0$& \cr
height .1cm &&height .1cm &&height .1cm &&height .1cm
  \cr
   \noalign{\hrule}
height .1cm &&height .1cm &&height .1cm &&height .1cm
  \cr
  &&&$b > e(a+1)$&&$0$& \cr
height .1cm &&height .1cm &&height .1cm &&height .1cm
  \cr    &&&\omit\hrulefill&&\omit\hrulefill&\cr
height .1cm &&height .1cm &&height .1cm &&height .1cm
  \cr
  &$a \leq -2$&&$b =e(a+1)$&&$?$&
 \cr
height .1cm &&height .1cm &&height .1cm &&height .1cm
  \cr &&&\omit\hrulefill&&\omit\hrulefill&\cr
height .1cm &&height .1cm &&height .1cm &&height .1cm
  \cr
  &&&$b < e(a+1)$&&$>0$&
  \cr
height .1cm &&height .1cm &&height .1cm &&height .1cm
  \cr}  \hrule }$$

\endproclaim

\proclaim {Proposition 5.3 (c.f. [GP] proposition 3.3)} The general
member of
$|C_0 -\frak e f|$ is a smooth elliptic curve and
those are the only smooth curves in the numerical
class of $C_0 +ef$.
\endproclaim

\vskip .2 cm
(5.3.1) We will fix once and for all a smooth elliptic curve $F$ in
the numerical class of $C_0 + ef$.
\vskip .1 cm

\proclaim {Proposition 5.4 ([H], V.2.21.b; [GP] Proposition 3.5 and
Remark 3.5.4)}
Let $L$ be line bundle on $X$ in the numerical class of $aC_0 +
bf$.
\roster
\item"5.4.1" $L$ is ample iff $a >0$ and $b >ae$
\item"5.4.2" $L$ is base point free if $a \geq 0$ and
$b-ae \geq 2$.
\item"5.4.3" $L$ is ample and base point free iff $a \geq 1$ and 	$b-ae \geq
2$.
\endroster
\endproclaim
We need another two lemmas in order to prove Proposition 5.1
\proclaim{Lemma 5.5}
Let $L =B_1 \otimes B_2$
and $L' = B_3$ be line bundles on $X$ satisfying the following
properties:
\roster
\item " 5.5.1" $B_1 \equiv B_3$
\item " 5.5.2" $B_i$ is in the numerical class of
$C_0 +(e+2)f$.
\endroster
Let $P$ and $P'$ be effective line bundles in the numerical classes of
$a(C_0+ef)+bf$ and of $a'(C_0+ef)+b'f$ respectively for some $a,b,a',b' \geq
0$. Then the map
$$\text H ^0 (L \otimes P) \otimes \text H^0(L' \otimes P')
@>\alpha >> \text H^0(L \otimes L' \otimes P \otimes P')$$
surjects.
\endproclaim

{\it Proof.} Analogous to the proof of Lemma 4.8. \boxit{}

\proclaim{Lemma 5.6}  Let $a$, $b$ be two integers such that $a \geq 1$ and
$b-ae \geq 4$. Let
$L$ be a line bundle in the numerical class of
$aC_0+bf$ and $P$ a line bundle in the numerical
class of $a'(C_0+ef)+b'f$ for some $a',b' \geq 0$. Then
$$\text H^2(M_L^{\otimes 2} \otimes P) =0$$ \endproclaim

{\it Proof.} Analogous to the proof of Lemma 4.9. \boxit{}

(5.7)
{\it Proof of Proposition 5.1.}

{\it Step 1.} $\text H^1 (M_{B^1} ^{\otimes p+1} \otimes
L_2)=0$

We will use Lemma 2.8.  The set
$\frak B$ will be the numerical class of $C_0 +(e+2)f$, the set $\frak P$
will consists of all line bundles numerically equivalent to $a(C_0+ef)+bf$
for some $a,b \geq 0$  and
$p_0$ will be equal to $1$. Therefore if $\frak B$ and $\frak P$ satisfy
the conditions of Lemma 2.8, we are done (simply take
$n$ to be $0$ and $P = P^2$  in the conclusion of the
Lemma 2.8). Conditions 2.8.1 and 2.8.2 are satisfied, (c.f.
Proposition 5.2 and 5.4). For Condition 2.8.3 note that $B_1 \otimes B_2
\otimes B_3^*$  belongs to $\frak P$ for any $B_1$,
$B_2$ and
$B_3$ in $\frak B$. In the other
hand, the line bundle $L_3$ defined in Lemma 2.8 satisfies the
hypothesis for the line bundle $L$ in Lemma 5.6. Thus applying
the mentioned lemma we are done. For Condition 2.8.4 we have to
show that
$$\text H^1(M_{B_1 \otimes B_2 \otimes C_1 \otimes
\dots \otimes C_n}^{\otimes 2} \otimes B_1' \otimes B_2'
\otimes P) = 0$$
for all $B_1, B_2, B_1', B_2', C_1, \dots ,C_n
\in
\frak B$ and $P \in \frak P$. Note that by 2.8.2 \ $C_1 \otimes
\dots \otimes C_n$ belongs to $\frak P$. Thus we will prove this more
general result:
$$\text H^1(M_{B_1 \otimes B_2 \otimes P^1}^{\otimes 2} \otimes
B_1' \otimes B_2'
\otimes P^2) = 0 \leqno (5.7.1)$$
for any $P^1$ and $P^2$ in $\frak P$.

First we show using Lemma 3.3 that
$$\text H^1(M_{B_1 \otimes
B_2 \otimes P^1}^{\otimes 2} \otimes B_1' \otimes B_2'
) = 0\ .\leqno (5.7.2)$$
We may assume without loss of generality that
$P^1 = \Cal O(a_1F+b_1f)$.
We use induction on $(a_1,b_1)$. If
$(a_1,b_1) = (0,0)$, the result follows from  Lemma 3.6.
Assume that
$$\text
H^1(M_{B_1 \otimes B_2 \otimes \Cal O((a_1-1)F)} \otimes
B_1' \otimes B_2') =0$$
We apply Lemma 3.3 to $B = B_1 \otimes B_2 \otimes \Cal
O((a_1-1)F)$, $P=\Cal O(C_0 +\frak e f)$, $\frak p = F$,
$L=B_1'
\otimes B_2'$, $p=1$, $p_0=-1$, $a=0$ and $m=0$.
Condition 3.3.1 is satisfied by Proposition 5.2. Condition 3.3.2 is
satisfied because deg$(B_1 \otimes B_2 \otimes \Cal
O_F(-lF)) \geq e +4 \geq 4$. We check that Condition 3.3.3 is
satisfied by using Lemma 4.7, noting that $$\text{deg}(B_1 \otimes
B_2
\otimes \Cal
O_F(a_1F)) > \text{deg}(B_1 \otimes B_2 \otimes \Cal
O_F((a_1-1)F) \geq (2 +a_1-1)e +4 \geq 4 \geq 4-l$$
and that
$$\text{deg}(B_1 \otimes
B_2 \otimes \Cal O_F(-lF)) \geq 4 \geq 3-l \ .$$
Condition 3.3.4 requires that $\text H^1(M_{B_1 \otimes B_2
\otimes \Cal O((a_1-1)F)} \otimes B_1' \otimes B_2' \otimes
\Cal O(-F))=0$, which is a consequence of Lemma 5.5,
and that
$$\text H^1(M_{B_1 \otimes B_2
\otimes \Cal O((a_1-1)F)}^{\otimes 2} \otimes B_1' \otimes B_2'
)=0$$
which is true by the induction hypothesis on $a_1-1$.
Condition 3.3.5 requires the vanishing of $\text H^1(B_1'
\otimes B_2' \otimes \Cal O(-2F))$ which follows from
Proposition 5.2.

To finish the proof of (5.7.2) we do induction on $b_1$. If
$b_1 =0$, the required statement has just been proven. Assume that
the result is true for $b_1 -1$ ($b_1 \geq 1$). We will use
again Lemma 3.3 setting $B = B_1^1 \otimes B_2^1
\otimes \Cal O(a_1(C_0-\frak e f)+(b_1-1)f)$, $L=B_1^2
\otimes B_2^2$, $P = \Cal O(f)$, $\frak p = f$,
$p=1$ and $p_0=-1$. Condition 3.3.1 is satisfied because
of Proposition 5.2. Condition 3.3.2 is satisfied because deg$(L_2
\otimes \Cal O_f(-lf)) \geq a_2 + 2\geq 2$.
Condition 3.3.3 follows from Lemma 5.6, since
$$\text{deg}(B_1^1
\otimes B_2^1
\otimes \Cal O_f(a_1(C_0+\frak e f)+(b_1-2)f)) \geq a_1 +
2 \geq 2
\qquad \text{and}$$
$$\text{deg}(L_2
\otimes \Cal O_f(-lf)) \geq a_2 + 2\geq 2 \ .$$
Condition 3.3.4 follows by induction hypothesis on $b_1-1$ and
from Lemma 5.5. Condition 3.3.5 follows
from Proposition 5.2.

To finish the proof of (5.7.1) we apply Lemma 3.2
inductively (as done for instance in the proof of Lemma 4.8) setting
the line bundles $B_1$ and
$B_2$ in the statement of Lemma 3.2 both equal to
$B_1
\otimes B_2
\otimes P^1$, $P$ equal
to $\Cal O(F)$ or $\Cal O(f)$,
$\frak p$ equal to $F$ or $f$ and $p = 2$.
\vskip .3 cm
{\it Step 2.} $\text H^1(M_{L_1} ^{\otimes p+1} \otimes L_2)
=0$.
We may assume without loss of generality that
 $L_1 = B_1^1
\otimes \dots \otimes B_{p+1}^1 \otimes \Cal O(a_1(C_0 -\frak
e f)+ b_1f)$. Thus we want to prove
$$\text H^1(M_{B_1^1
\otimes \dots \otimes B_{p+1}^1 \otimes \Cal O(a_1(C_0 -\frak
e f)+ b_1f)} ^{\otimes p+1} \otimes L_2)
=0\ .\leqno (5.7.3)$$
We will use
induction on $p$, starting at $p=1$. If
$p=1$ (5.7.3) follows from (5.7.1).

Now we assume that (5.7.3) holds for $1, \dots ,p-1$ for $p
\geq 2$ and we will prove that it holds also for $p$. Again we do
induction on $(a_1, b_1)$. If $(a_1,b_1) = (0,0)$ the
statement was proven in Step 1. Assume the result is true for
$(a_1-1,0)$. We apply Lemma 3.3 to $B=B^1 \otimes \Cal
O((a_1-1)(C_0 + \frak e f))$, $L = L_2$, $\frak p = F$ and
$p_0=-1$. Condition 3.3.1 is satisfied by Proposition 5.2. Condition
3.3.2 follows from the fact that deg$(L_2 \otimes \Cal
O_{F}(-lF)) \geq (p+1)(e+2) + (a_1-1)e \geq 2p+2 >0$.
Condition 3.3.3 follows from Lemma 4.7 and from the fact that
deg$(L_2
\otimes
\Cal O_{F}(-lF)) \geq 2p+2 > p+2$ and $$\text{deg}(B^1
\otimes
\Cal O_F(aF) \geq \text{deg}(B^1 \otimes \Cal
O_F((a_1-1)F) \geq 2p+2 \geq p+3 \ .$$
Condition 3.3.4 requires the vanishing of
$$\text H^1(M_{B^1
\otimes \Cal
O((a_1-1)F)} ^{p'+1} \otimes L_2 \otimes \Cal O(-lF))
\leqno (5.7.4)$$ for
$l=p-p'$ and $0 \leq p' \leq p$. If $p = p'$, (5.7.4) is
simply the induction hypothesis for $a_1-1$. If $1 \leq p'
\leq p-1$, (5.7.4) is nothing but the induction hypothesis on
$1, \dots , p-1$. Indeed.
The line
bundle
$L_2
\otimes
\Cal O(-lF)$ can be written as the tensor product of $B_1^2
\otimes \dots
\otimes B_{p'+1}^2$ with an effective line bundle in the
numerical class of $a_2(C_0+ef)+(b_2 + 2(p-p'))f$. If $p' =0$,
(5.7.4) follows from Lemma 5.5. Condition 3.3.5
requires the vanishing of
$\text H^1(L_2 \otimes \Cal O(-(p+1)F))$ which follows from
Proposition 5.2.

The induction argument on $b_1$ is similar to the one on
$a_1$. \boxit{}

\heading 6. Syzygies of elliptic ruled surfaces \endheading

In this section we assume that $\text{char}(\bold k) >p+1$ or
equal to $0$. We will use the results obtained in Sections 4 and 5 to prove the
following

\proclaim{Theorem 6.1 }
Let
$X$ be an elliptic ruled surface and  let $p \geq 1$. Let $a$, $b$
be integers and let
$L$ be a line bundle in the numerical class of $aC_0 + bf$.
\roster
\item "6.1.1"
If
 $e = e(X) = -1$ and
$a \geq p+1$, $a+b \geq 2p+2$ and $a+2b \geq 2p+2$, then
$L$ satisfies the property $N_p$.
\item "6.1.2"
If $e =e(X) \geq 0$ and $a \geq p+1$, $b-ae \geq 2p+2$, then
$L$ satisfies the property $N_p$.
\endroster
\endproclaim

(6.1.3) Note that if $p=1$, we recover from Theorem 6.1, the ``if"
part of Theorem 4.2 of [GP], except for the case when $a=1$.
\vskip .4 cm
(6.1.4) {\it Proof of Theorem 6.1.} The line bundle $L$ is normally
generated (see [Ho1] and [Ho2]; see also [GP], Lemma 2.6 and Theorem 4.2).
Hence by Lemma 1.2 and 1.2.1 (this  is the reason why we need the
hypothesis on the characteristic of $\bold k$),
it is enough to show  that
$$\text H^1(M_L^{\otimes k+1} \otimes L) =0 \ \text{for all}\
1
\leq k
\leq p \ .\leqno (6.1.5)$$

If $e =-1$, $L$ can be written for all $1 \leq k \leq p$
either as $B_1 \otimes \dots \otimes B_{k+1} \otimes P$, where
$B_i$ is in the numerical class of $2C_0$ or of $C_0+f$ and $P$
is effective in the numerical class of $aC_0+bf$ for some $a,b
\geq 0$ or as $B_1 \otimes \dots \otimes B_{k+1} \otimes P$, where
$B_i$ is in the numerical class of $2C_0$ and $P$
is effective. Thus by Proposition 4.1 and Proposition 4.2 we
obtain the result.

If $e \geq 0$, $L$ can be written for all $1 \leq k \leq p$  as $B_1 \otimes
\dots \otimes B_{k+1} \otimes P$, where
$B_i$ is the numerical class of $C_0 +(e+2)f$ and $P$
is effective in the numerical class of $a(C_0+ef)+bf$ for some
$a,b
\geq 0$. Thus by Proposition 5.1  we
obtain the result. \boxit{}
\vskip .2 cm
As a corollary of Theorem 6.1 we obtain the following result on
adjoint linear series, which is a generalization to higher
syzygies of Corollary 4.6 of [GP]. Note however that we obtain there a
sharper bound in the case $e\geq 1, p=1$.

\proclaim{Corollary 6.2} Let $X$ be an elliptic
ruled surface and let $p \geq 1$. Let $A_i$ be  an
ample line bundle on $X$ for all $1 \leq i \leq q$.
If $q \geq 2p+2 - \text{min}(e(X),p-1)$, then
$\omega_X \otimes A_1 \otimes \dots \otimes A_q$
satisfies the property $N_p$. \endproclaim

{\it Proof.}
Let $A_i$ be in the numerical class of $a_iC_0 +b_if$
and $\omega _X \otimes A_1 \otimes \dots \otimes A_q$
in the numerical class of $aC_0+bf$. If $e=-1$, $A_i$ is ample
iff $a_i \geq 1$ and $a_i +2b_i \geq 1$ (c.f. Proposition 4.5). In
particular we also have that if $A_i$ is ample, then
$a_i+b_i \geq 1$. Since $\omega _X $ is numerically
equivalent to $-2C_0+f$ it follows that
$$\displaylines {a\geq q-2 \geq 2p+1 >p+1 \cr
a+b \geq q-1 \geq 2p+2 \quad \text{and} \cr
a+2b \geq q \geq 2p+3 \ .}$$
Hence by Theorem 6.1, $\omega _X \otimes A_1 \otimes \dots
\otimes A_q$ satisfies the property $N_p$.

If $e \geq 0$ $A_i$ is ample
iff $a \geq 1$ and $b_i - a_ie \geq 1$ (c.f. Proposition 5.4).
Since
$\omega _X $ is numerically
equivalent to $-2C_0-ef$ it follows that $a \geq q-2$
and $b-ae \geq q+e$. By hypothesis, $q \geq 2p+2-e$ and
$q \geq p+3$; hence $a \geq p+1$ and $b-ae \geq 2p+2$
and by Theorem 6.1, $\omega _X \otimes A_1 \otimes \dots
\otimes A_q$ satisfies the property $N_p$. \boxit{}
\vskip .2 cm
We also obtain this generalization of Corollary 4.4 of [GP]:

\proclaim{Corollary 6.3} Let $X$ be as above and let $p
\geq 1$. Let $B_i$ be an ample and base-point-free
line bundle on $X$ for all $1 \leq i \leq q$. 	If
$q \geq p+1$, then $B_1 \otimes \dots \otimes B_q$
satisfies the property $N_p$. \endproclaim

{\it Proof.} Let $B_i$ be in the numerical class of $a_iC_0 +
b_if$ and $B_1 \otimes \dots \otimes B_q$ in the numerical
class of $aC_0+bf$. If $e=-1$, by Proposition 4.5,
$B_i$ is ample and base-point-free iff $a_i \geq 1$, $a_i +
b_i \geq 2$ and
$a_i+2b_i
\geq 2$. Thus we obtain that
$$\displaylines{a \geq q \geq p+1 \cr
a+b \geq 2q \geq 2p+2 \cr
a+2b \geq 2q \geq 2p+2 }$$
Hence, by Theorem 6.1, $B_1 \otimes \dots \otimes B_q$
satisfies the property $N_p$.

If $e \geq 0$, $B_i$ is ample and base-point-free iff $a_i
\geq 1$ and $b_i -a_ie \geq 2$ (c.f. Proposition 5.4). Thus $a
\geq q \geq p+1$ and $b-ae \geq 2q \geq 2p+2$. Hence, by
Theorem 6.1, $B_1 \otimes \dots \otimes B_q$ satisfies the
property $N_p$.
\boxit{}

\proclaim{Corollary 6.4} Let $X$ as above and let $p
\geq 1$. Let $A_i$ be an ample
line bundle on $X$ for all $1 \leq i \leq q$. 	If
$q \geq 2p+2$, then $A_1 \otimes \dots \otimes A_q$
satisfies the property $N_p$.\endproclaim

{\it Proof.} If suffices to note that if $A$ and $A'$ are ample line bundles on
$X$, then $A \otimes A'$ is ample and base-point-free (this follows
from Propositions 4.5 and 5.4). Then we apply Corollary 6.3.
\boxit{}

(6.5) \ Note that the assumption on the characteristic was made
because we wanted to be able to consider
$\bigwedge^{p'+1} M_L \otimes L^{\otimes p'+1}$ as a direct summand
of
$M_L^{\otimes p'+1}\otimes L^{\otimes p'+1}$, for all $1 \leq p'
\leq p$. That way we obtained from the vanishings of $\text
H^1(M_L^{\otimes p'+1}\otimes L^{\otimes p'+1})$, the vanishings
of $\text H^1(\bigwedge^{p'+1} M_L \otimes L^{\otimes p'+1})$, for
all $1 \leq p'
\leq p$. These were the vanishings required by Lemma 1.2 in order
that
$L$ satisfied the property $N_p$. However, in particular situations,
those conditions required in Lemma 1.2 can be relaxed. Precisely,
if $L$ is a normally generated line bundle such that $\text
H^i(L^{\otimes 2-i})=0$ and $p$ is less or equal than the
codimension of $X$ inside $\bold P^N = \bold P(\text H^0(L))$, then
$L$ satisfies the property $N_p$ iff the group $\text
H^1(\bigwedge^{p+1} M_L
\otimes L^{\otimes p+1})$ vanishes (c.f [GL], Lemma 1.10). We
claim that the above condition on $p$ and the codimension is
satisfied under the conditions of Theorem 6.1.

Indeed. If $L$ belongs to the numerical class of $aC_0+bf$, using
Riemann-Roch one easily obtains that
$\text h^0(L)= \frac 1 2(a(b-1)+(a+2)-a(a+2)e)$. Thus, if $e
=-1$, we want to see that $$\frac{(a+1)(a+2b)}{2}-3
= \text{cod}(X,
\bold P^N)
\geq p \ .$$ The latter inequality follows from the numerical
conditions satisfied by $(a,b)$:
$$\frac{(a+1)(a+2b)}{2}-3\geq (p+2)(p+1)-3\geq p\ ,\quad
\text{for all}\ p \geq 1.$$
If, on the other hand, $e \geq 0$, using again the numerical
conditions satisfied by $(a,b)$, we see that
$$\displaylines{\hskip .2 cm\text{cod}(X,
\bold P^N)=\frac{(a+2)(b-ae)+ a(b-1)} 2 -3 \hfill\cr\hfill
\geq
\frac{(p+3)(2p+2)+3(p+1)} 2 -3 \geq p \hskip .1 cm}$$ for all $p
\geq 1$.
\vskip .1 cm
Hence the
results of this section hold in slightly greater generality,
namely, they hold when
$\text{char}(\bold k)$ does not divide $p+1$.

\heading 7. Open questions and conjectures \endheading
We foresee two directions in which these results on
 syzygies of elliptic ruled surfaces could be improved:
\vskip .15 cm
{\bf 7.1.} In Section 4 of [GP] we prove that the product of two
base-point-free divisors (not necessarily both of them ample)
satisfies the property $N_1$ iff it is ample. Therefore one may ask whether a
similar statement is true for any $p \geq 1$, i.e., whether the product $L$ of
$p+1$ base-point-free divisors (not necessarily all of them ample) satisfies
the
property $N_p$ whenever $L$ is ample. This is expressed graphically for the
case
$e(X)=-1$ in Figure 1 and  for the case $e \geq 0$ in Figure 2. In these
figures the integral points of the coordinate plane represent the classes of
$\text{Num}(X)$) and the  shadowed regions contain the divisors
which could satisfy the property
$N_p$.
\vskip .15 cm
(7.2) When $e(X)=-1$, Homma proved (see [Ho2]) that
if $L$ is a line bundle in the numerical class of $aC_0+bf$, then $L$ satisfies
the property $N_0$  iff $a \geq 1$, $a+b \geq 3$, and $a+2b \geq 3$. We prove
in Theorem 4.2 of [GP] that $L$ satisfies the property $N_1$ iff  $a \geq 1$,
$a+b \geq 4$, and $a+2b \geq 4$. Hence one could ask whether $L$ satisfies the
property $N_2$ if $a \geq 1$,
$a+b \geq 5$, and $a+2b \geq 5$. Evidence suggesting an affirmative answer is
the fact that the free resolution of $R(L)$ is
linear until the second stage if $L$ is in the numerical class of $5f$ and
if $L$ is certain line bundle in the numerical class of $C_0+4f$ and in the
class of $2C_0+3f$ (these two cases were checked  using the computer program
Macaulay). Analogously, one expects similar statements for $p \geq 3$ and
also for the case $e \geq 0$. We make the following

\proclaim {Conjecture 7.3} Let $X$ be an elliptic ruled surface and let $L$
be a line bundle on $X$ in the numerical class $aC_0+bf$.

If
$e=-1$,
$L$ satisfies the property
$N_p$ iff $a \geq 1$,
$a+b \geq p+3$, and $a+2b \geq p+3$.

If $e(X) \geq 0$,  $L$ satisfies the
property $N_p$ iff $a \geq 1$ and $b-ae \geq p+3$.
\endproclaim

In Figure 3 we show, for the case $e=-1$ the lines (dashed) joining the
numerical classes of those line bundles which are conjectured to be optimal
``$N_p$ line bundles". If  this conjecture is true,
$\omega_X \otimes A^{\otimes p+4}$ will
satisfy the property $N_p$.
Hence Conjecture 7.3 implies Mukai's conjecture in the case of elliptic ruled
surfaces. It also implies an affirmative answer for Question 7.1.
\vskip .2 cm
(7.4) Observe the analogy of Conjecture 7.3 and Green's Theorem for curves,
which says that $L$ satisfies the property $N_p$ if $\text{deg} (L) \geq
2g+p+1$. There the difference between two  consecutive bounds is $1$, i.e., the
minimal degree for an ample line bundle on a curve. Going back to elliptic
ruled surfaces, the ``difference" between the line joining the conjectured
optimal ``$N_p$ line bundles" and the line joining the conjectured optimal
``$N_{p+1}$ line bundles" is $C_0$, which is the ``minimal" ample divisor.

\vfill\eject
\epsfbox{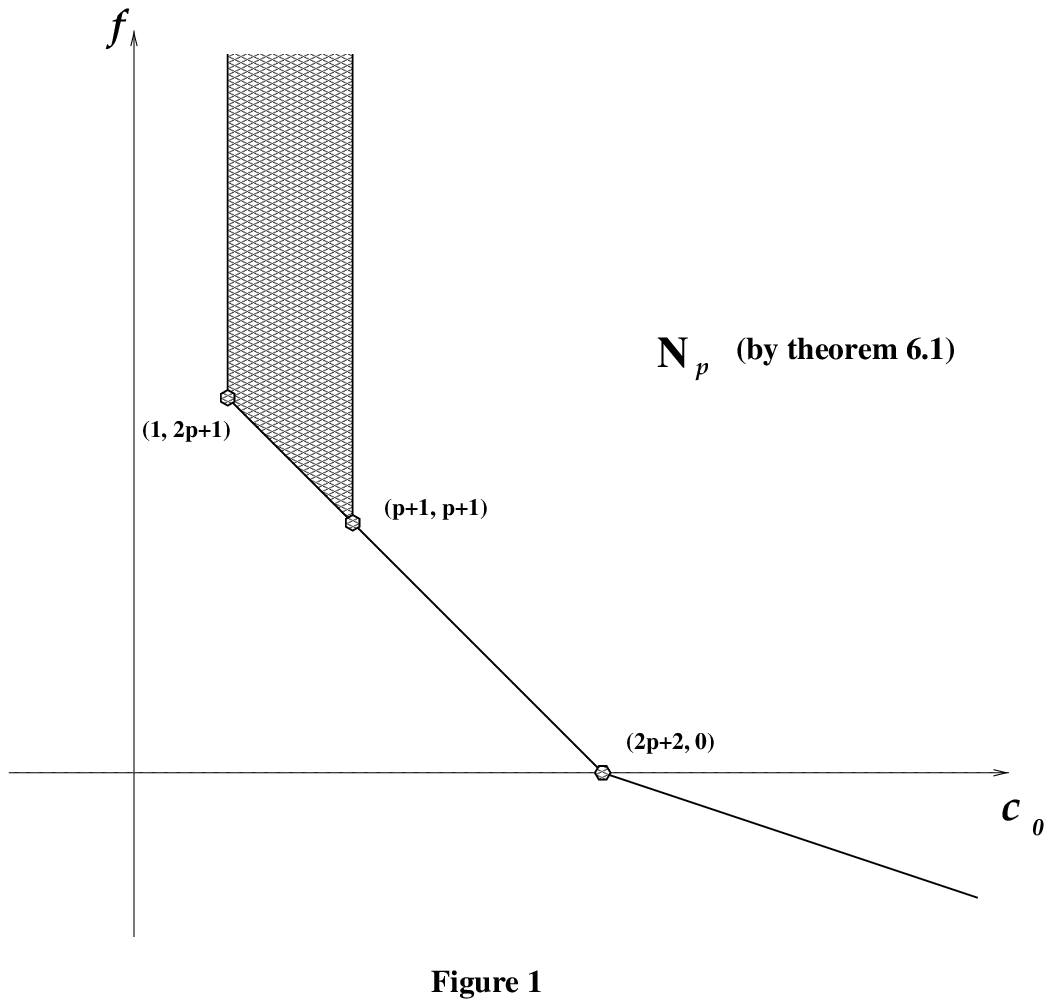}
\epsfbox{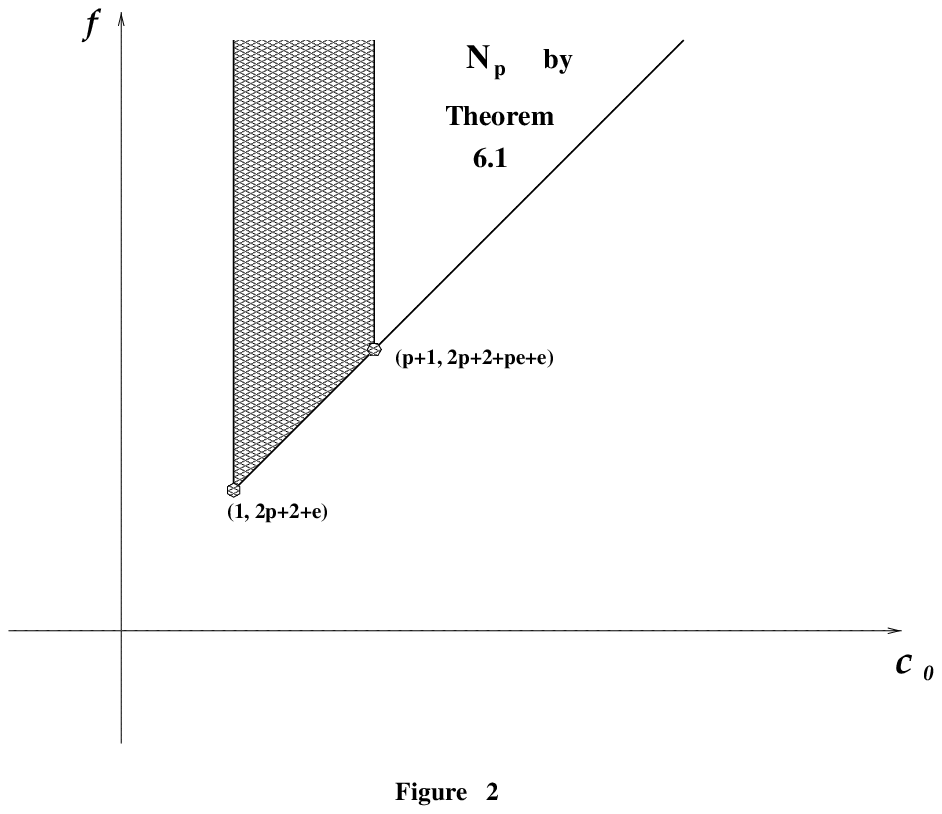}
\vfill\eject
\epsfbox{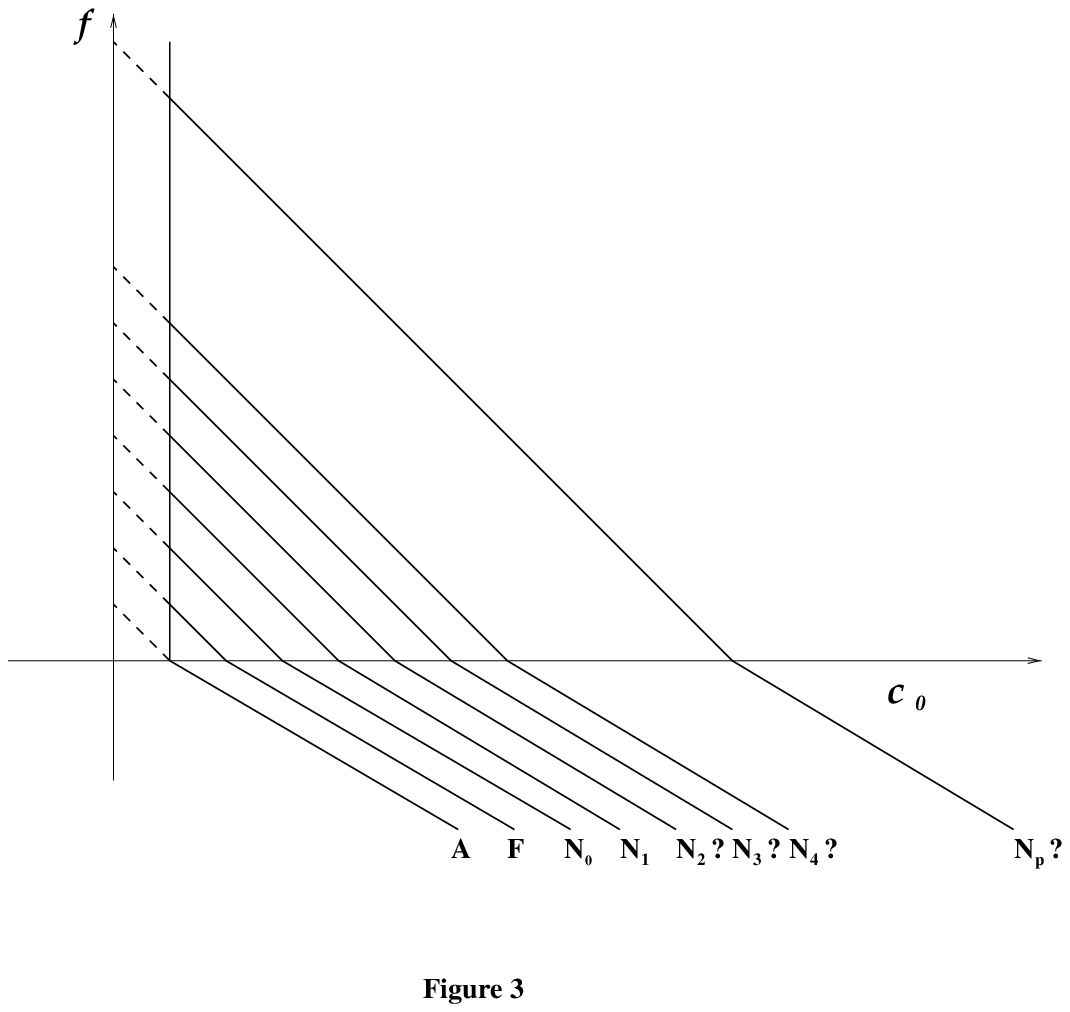}
\heading References  \endheading
\roster

\item
"[B]" E. Bombieri, {Canonical models of surfaces of general type}, IHES,
{\bf 42} (1973)

\item
"[Bu]" D. Butler, {\it Normal generation of vector bundles over
a curve}, J. Differential
Geometry {\bf 39} (1994) 1-34.

\item"[G]" M. Green, {\it Koszul cohomology and the geometry of
projective varieties},
J. Differential Geometry {\bf
19} (1984) 125-171.

\item"[GL]" M. Green \& R. Lazarsfeld, {\it Some results on the
syzygies of finite sets and
algebraic curves}, Compositio
Math. {\bf 67} (1989) 301-314.

\item "[GP]" F.J. Gallego \& B.P. Purnaprajna, {\it Normal
Presentation on Elliptic Ruled Surfaces} (1995) Preprint.

\item "[GP1]" F.J. Gallego \& B.P. Purnaprajna, {\it Syzygies
of K3 surfaces and Fano varieties} (1995) Preprint.

\item "[GP2]" F.J. Gallego \& B.P. Purnaprajna, {\it Syzygies
of surfaces and Calabi-Yau threefolds} (1995) In preparation.

\item"[H]" R. Hartshorne {\it Algebraic Geometry}, Springer,
Berlin, 1977.

\item"[Ho1]" Y. Homma, {\it Projective normality and the
defining equations of ample
invertible sheaves on elliptic ruled
surfaces with $e \geq 0$}, Natural Science
Report, Ochanomizu
Univ. {\bf 31} (1980) 61-73.

\item"[Ho2]" \hbox{\leaders \hrule  \hskip .6 cm}\hskip .05 cm , {\it
Projective normality and the defining equations of an elliptic ruled
surface with negative
invariant}, Natural Science Report, Ochanomizu
Univ. {\bf 33} (1982) 17-26.

\item"[L]" R. Lazarsfeld, {\it A sampling of vector bundles techniques in the
study of linear series}, Lectures on Riemann Surfaces, World Scientific
Press, Singapore, 1989, 500-559.

\item"[Mi]" Y. Miyaoka, {\it The Chern class and Kodaira
dimension of a minimal variety},
Algebraic Geometry --Sendai 1985,
Advanced Studies in Pure Math.,
Vol. 10,
North-Holland,
Amsterdam, 449-476.

\item "[Mu]" D. Mumford, {\it Varieties defined by quadratic
equations}, Corso CIME in
Questions on Algebraic Varieties,
Rome, 1970, 30-100.

\item"[R]" I. Reider, {\it Vector bundles of rakk $2$  and
linear systems on an algebraic surface},
Ann. of Math. (2) {\bf
127} (1988) 309-316.

\endroster

\enddocument